\documentclass{amsart}
\usepackage{amssymb}
\input xypic
\xyoption{all}

\newtheorem{theorem}{Theorem}[section]

\newtheorem{proposition}[theorem]{Proposition}

\newtheorem{lemma}[theorem]{Lemma}
\newtheorem{corollary}[theorem]{Corollary}
\theoremstyle{definition}
\newtheorem{definition}[theorem]{Definition}
\newtheorem{remark}[theorem]{Remark}

\numberwithin{equation}{section}
\newcommand{\mer}{\mathcal{M}}
\newcommand{\T}{\mathcal{T}}
\newcommand{\Pic}{\operatorname{Pic}}
\newcommand{\divar}{\operatorname{Div_{Ar}}}
\newcommand{\chowar}{\operatorname{CH_{Ar}}}
\newcommand{\picar}{\operatorname{Pic_{Ar}}}
\newcommand{\prinar}{\operatorname{Prin_{Ar}}}
\newcommand{\Spec}{\operatorname{Spec}}
\newcommand{\Gal}{\operatorname{Gal}}
\renewcommand{\det}{\operatorname{det}}
\newcommand{\dd}{\textup{d}}
\newcommand{\qpc}{\overline{\Q_p}}
\newcommand{\fvp}{\overline{F_v}}
\newcommand{\alog}{A_{\textup{log}}}
\newcommand{\Olog}{\Omega_{\textup{log}}}
\newcommand{\aloga}{A_{\textup{log},1}}
\newcommand{\Ker}{\operatorname{Ker}}
\newcommand{\C}{{\mathbb{C}}}
\newcommand{\R}{{\mathbb{R}}}
\newcommand{\NN}{\mathcal{N}}
\newcommand{\PP}{{\mathbb{P}}}
\newcommand{\ninv}{\Psi}
\newcommand{\inject}{\hookrightarrow} 

\newcommand{\Ext}{\operatorname{Ext}}
\newcommand{\res}{\operatorname{Res}}
\newcommand{\norm}{\operatorname{Norm}}
\newcommand{\xvb}{\overline{X_v}}
\newcommand{\symm}{\operatorname{Symm}}
\newcommand{\id}{\operatorname{Id}}
\newcommand{\dR}{{\textup{dR}}}
\newcommand{\hdr}{H_{\dR}}
\newcommand{\dlog}{\operatorname{d\!\log}}
\newcommand{\dllog}{\operatorname{d\!\llog}}
\newcommand{\ilog}{\operatorname{\iota_{\log}}}
\newcommand{\pair}[1]{{\left\langle #1 \right\rangle}}
\newcommand{\nekovar}{Nekov\'a\v r}
\renewcommand{\O}{\mathcal{O}}
\newcommand{\M}{\mathcal{M}}
\newcommand{\pr}{^\prime}
\newcommand{\Div}{\operatorname{div}}
\newcommand{\lag}{\Psi}
\newcommand{\indi}{\rho}
\newcommand{\bfo}{{\boldsymbol{\omega}}}
\newcommand{\gl}{{\textup{gl}}}
\newcommand{\db}{\bar{\partial}}
\newcommand{\dbd}{\db\partial}
\newcommand{\Lie}{\operatorname{Lie}}
\newcommand{\tot}{{\operatorname{Tot}}}
\renewcommand{\L}{\mathcal{L}}
\newcommand{\Curve}{\operatorname{Curve}}
\newcommand{\Ht}{H^\otimes}
\def\htt_#1{H_#1^\otimes}
\newcommand{\Lt}{\L^\times}
\newcommand{\llog}{\log_{\L}}
\newcommand{\col}{\textup{Col}}
\newcommand{\ocol}{\Omega_{\col}}
\newcommand{\Hom}{\operatorname{Hom}}
\newcommand{\acol}{\O_{\col}}
\def\acoln_#1{\O_{\col,#1}}
\def\Ocoln_#1{\Omega_{\col,#1}^1}
\newcommand{\shF}{\mathcal{F}}
\newcommand{\shG}{\mathcal{G}}
\newcommand{\XX}{\mathcal{X}}
\newcommand{\im}{\operatorname{Im}}
\newcommand{\tr}{\operatorname{tr}}
\newcommand{\cech}{\v Cech}
\newcommand{\bo}{{\bar{\omega}}}
\newcommand{\isom}{\cong}
\newcommand{\Q}{\mathbb{Q}}
\newcommand{\ndiv}{\nmid}
\newcommand{\fin}{{\textup{fin}}}
\def\DD^#1{\Delta^{[#1]}}
\def\iny(#1,#2){{#2}^{-\infty}#1}
\newcommand{\tknd}[1]{\mathcal{C}_{#1}}
\newcommand{\sknd}[1]{\mathcal{B}_{#1}}
\begin{document}
\title[$p$-adic Arakelov theory]{$p$-adic Arakelov theory}
\author{Amnon Besser}
\address{Department of Mathematics\\
Ben-Gurion University of the Negev\\
P.O.B. 653\\
Be'er-Sheva 84105\\
Israel
}
\email{bessera@math.bgu.ac.il}
\subjclass{Primary 14G40; Secondary 11S80} 
\keywords{Arakelov theory, $p$-adic height pairings, $p$-adic
integration, $p$-adic Green functions}

\begin{abstract}
We introduce the $p$-adic analogue of Arakelov intersection
theory on arithmetic surfaces. The intersection pairing in an
extension of the $p$-adic
height pairing for divisors of degree $0$ in the form described by Coleman
and Gross. It also uses Coleman integration and is related to work of
Colmez on $p$-adic Green functions.  We introduce the $p$-adic version of
a metrized line bundle and define the metric on the determinant of its
cohomology in the style of Faltings. It is possible to prove in this
theory analogues of the Adjunction formula and the Riemann-Roch formula.
\end{abstract}

\maketitle

\section{Introduction}
\label{sec:introduction}

The purpose of this paper is to create a $p$-adic analogue to the part
of Arakelov theory that deals with arithmetic surfaces~\cite{Arak74,Fal84}.

In the classical case of Arakelov theory the
usual intersection pairings above finite primes are supplemented by a
pairing at infinity, involving analysis on the resulting Riemann
surface. Likewise, in the $p$-adic theory the same pairings at primes
not above $p$ are supplemented by a pairing at $p$ involving ``Coleman
analysis''~\cite{Col-de88,Colm96,Bes99} on the completions above
$p$. This ``local'' part of the
theory pauses most of the difficulties while using it to produce
global results is fairly straightforward given the classical case.

We would like to clarify that this work has nothing to do with the
work of Bloch, Gillet and Soule on ``non-archimedean Arakelov
theory''. Much like with the term ``$p$-adic integration'', which is
confusingly used for a number of non-related constructions, there is
ample room for confusion here. The construction of non-archimedean
Arakelov theory give real valued results, whereas $p$-adic Arakelov
theory gives $p$-adic results.

Our starting point was the theory of $p$-adic height pairings,
especially in the form of Coleman and Gross in~\cite{Col-Gro89}. The local part
at $p$ of that theory gives the intersection index $\pair{D,E}$ for
two divisors $D$ and $E$, of degree $0$, on a complete curve over a
$p$-adic field. The goal of the local theory is to describe an
extension of this pairing, by giving a Green function $G(P,Q)$
with certain desirable properties, in such a way that
\begin{equation}\label{heightpairing}
  \pair{\sum n_i P_i,\sum m_j Q_j}=\sum n_i m_j G(P_i,Q_j)\;.
\end{equation}

To isolate a canonical as possible Green function, one needs some
extra conditions. The most natural comes by introducing a notion of
metrized line bundles. These are line bundles together with a function
which behaves like a log (up to scaling) on the fibers and which is a
Coleman function of a certain type. Having this notion one can impose
the analogue of the condition, which is satisfied by the canonical
Green function in the classical theory, that the residue map defines a
metric on the canonical bundle, which is admissible with respect to
the Green function (or rather, the volume form). This extra condition
indeed isolates a canonical choice of a Green function, up to a
constant \eqref{Green-formula}. It is possible to work with this
definition but this becomes very cumbersome. We have therefore chosen
another route, which runs in parallel with the classical theory.

We utilize the definition of the $p$-adic $\db$ operator
from~\cite{Bes99} to define the curvature of a metrized line
bundle. Once this is done, the Green function on a curve $X$ is
derived from the metric on the line bundle $\O(\Delta)$, where
$\Delta$ is the diagonal in $X\times X$, having a prescribe curvature
similar to the one encountered in the classical theory. Unlike the
classical theory we still need to impose the residue condition
described before to obtain a unique choice up to constant. The
relation with the Coleman-Gross height pairing now requires a proof.

The advantage of this approach is that we have a better understanding
of the Green function as a function of two variables. This allows us
to define the analogue of the Faltings volume on the determinant of
cohomology.

A $p$-adic Green function on curves was previously defined by Colmez
in~\cite{Colm96}, using abelian varieties. We related his theory with
ours in Appendix~\ref{sec:colmez}.

We then turn to the global theory, proving analogues of the adjunction
formula and the Riemann-Roch formula. We are unfortunately unable at
the moment to produce a result corresponding to the Noether
formula because of the problem with the normalization of the Green
function.

Still missing are applications. Unlike real heights, it is not clear
to us what $p$-adic heights are good for except exact formulas. It
seems to us though that research in classical Arakelov theory is leaning more
and more towards exact formulas involving heights. Such results could
eventually find $p$-adic analogues.

This word begun when the author was at the university of M\"unster and
some crucial progress was made while visiting IHES\@. We would like to
thank both institutions.
We would also like to thank Amnon Yekutieli.

\section{Review of $p$-adic integration}
\label{sec:review}

The $p$-adic analysis required for defining the part of Arakelov
intersection ``at infinity'', i.e., at $p$, is given by the theory of
Coleman integration. We use here the tools developed
in~\cite{Bes99}. However, the work of Vologodski~\cite{Vol01}, is
extremely useful here because it applies to the case of bad reduction
as well, and because it works with the Zariski topology instead
of with the rigid topology, which is particularly convenient for
working with Arakelov geometry (one drawback is that it works only
over finite extensions of $\Q_p$, but this is all we need). Our goal
in this section is to recall
the constructions of~\cite{Bes99} while showing how they work in the
context of Vologodski's work. The proofs are mostly easy modifications
of the ones in~\cite{Bes99} so we do not repeat them here.

We let $K$ be a finite extension of $\Q_p$. We
choose a branch of the logarithm. In~\cite{Vol01} a branch is not
chosen and the integration takes place in a ring containing a formal
variable for the log of $p$. In our setup we will choose a branch
anyhow and the specialization to this situation is clear.

Let $X$ be a smooth, geometrically connected algebraic variety over
$K$. Let $(F,\nabla)$ be a
unipotent connection over $X$, i.e., a coherent sheaf $F$ with an
integrable connection $\nabla$ which is a successive extension of
trivial connections. The main result of Vologodski is
\begin{theorem}[{\cite[Theorem B]{Vol01}}]\label{volsthm}
  for any two points $x$, $y\in X(K)$ there exists a canonical parallel
  translation isomorphism $v_{x,y}=v_{x,y}^F:F_x \to F_y$ of the fibers of $F$
  over $x$ and $y$. This translation satisfying the following properties:
  \begin{enumerate}
  \item The translation $v_{x,x}$ is the identity.
  \item For any 3 points $x,y,z$ we have $v_{y,z}\circ
    v_{x,y}=v_{x,z}$.
  \item \label{vol3} The translation $v_{x,y}$ is locally analytic in
    $x$ and $y$.
  \item For the trivial connection on $\O_X$ the translation carries
    $1$ to $1$.
  \item \label{vol5} For any map $T:E\to F$ of unipotent connections
    let $T_x$ and
    $T_y$ be the restrictions of $T$ to the fibers at $x$ and
    $y$. Then we have $T_y\circ v_{x,y}^E = v_{x,y}^F \circ T_x$.
  \item For any two unipotent connections $E$ and $F$ we have
    $v_{x,y}^{E\otimes F}=v_{x,y}^E \otimes v_{x,y}^F$.
  \item For any $K$-morphism $f: X'\to X$, and unipotent connection
    $F$ on $X$ and any two points $x,y\in X'(K)$ we have 
    $v_{x,y}^{f^\ast F} \circ (f^\ast)_x = (f^\ast)_y \circ
    v_{f(x),f(y)}^F$, where $(f^\ast)_x: F_{f(x)} \to (f^\ast F)_x$ is
    the pullback map (similarly with $x$ replaced by $y$).
  \item \label{vol8} Let $\sigma:\Spec(L) \to \Spec(K)$ be a finite
    map, where $L$
    is another field. Let $X_L:= X\times_{\Spec(K)} \Spec(L)$ and let,
    for a unipotent connection $F$ on $X$, $F_L$ denote the extension
    of scalars. For $x\in X(K)$ let $\sigma(x)\in X_L(L)$ be the
    corresponding point and let $\sigma:F_x \to (F_L)_{\sigma(x)}$ be
    the obvious map. Then the parallel translation is compatible with
    $\sigma$ in the sense that $\sigma \circ v_{x,y}^F =
    v_{\sigma(x),\sigma(y)}^{F_L} \circ \sigma$.
  \end{enumerate}
\end{theorem}
Note that property \eqref{vol5} is not stated as such in loc.\ cit.\ but it
follows from other properties. Also note that the locally analytic
nature of $v_{x,y}$ means that if $s\in F_x$, then $r(y)=v_{x,y} s$ is
a locally analytic section of $F$.

\begin{definition}[{Compare~\cite[Definition~4.1]{Bes99}}]
  A (Vologodski) \emph{abstract Coleman function} on $X$ with values in a
  locally free sheaf $\shF$ is a fourtuple
  $(M,\nabla,s,y)$ consisting of a unipotent connection $(M,\nabla)$
  on $X$ together with a homomorphism $s\in \Hom(M,\shF)$ (a sheaf but
  not a connection homomorphism) and a compatible system $y=(y_x)$ of
  elements $y_x\in M_x$ for every $x\in X(L)$, for any finite
  extension $L$ of $K$. This system should satisfy the following two
  conditions
  \begin{enumerate}
  \item For any two points $x_1,x_2\in X(L)=X_L(L)$, parallel translation on
    $X_L$ takes $y_{x_1}$ to $y_{x_2}$.
  \item For any map of fields $\sigma:L_1 \to L_2$ fixing $K$ for any
    point $x\in X(L_1)$ we have $\sigma(y_x)=y_{\sigma(x)}$.
  \end{enumerate}
  A morphism between two abstract Coleman
  functions with values in $\shF$, $(M_i,\nabla_i,s_i,y_i)$, $i=1,2$, is a map
  $f:(M_1,\nabla_1)\to (M_2,\nabla_2)$ pulling back $s_2$ to $s_1$ and
  sending $y_1$ to $y_2$. A \emph{Coleman function} with values in $\shF$ is
  a connected component of the category of abstract Coleman
  functions. We denote the connected component of $(M,\nabla,s,y)$ by
  $[M,\nabla,s,y]$. The collection of Coleman functions on $X$ with
  values in $\shF$ is denoted by $\acol(X,\shF)$. In particular we
  have set $\acol(X):=\acol(X,\O_X)$ and
  $\ocol^i(X):=\acol(X,\Omega_{X/K}^i)$.
\end{definition}

A Coleman function can be interpreted as a set theoretic section of
$\shF$ over $X(\bar{K})$ as follows: if $f$ corresponds to
$(M,\nabla,s,y)$ then $f(x)=s(y_x)$. From the definition of Coleman
functions it is clear that this function depends only on the
underlying Coleman function and not on the abstract Coleman function
used to define it. It follows from part~\ref{vol3} of
Theorem~\ref{volsthm} that Coleman functions are
locally analytic. From part~\ref{vol8} it follows that Coleman
functions are Galois
equivariant over $K$ in the 
sense that $f(\sigma(x))=\sigma(f(x))$ for any automorphism of fields
fixing $K$. There is a structure of a
$K$-vector space on $\acol(X,\shF)$ and a bilinear product
$\acol(X,\shF)\times \acol(X,\shG)\to \acol(X,\shF\otimes \shG)$ given
by obvious operations on abstract Coleman functions
(compare~\cite[Definition~4.2]{Bes99}). These operation correspond to
operations of
addition and multiplication when computed on points.

There are differentials $\dd: \ocol^i(X) \to \ocol^{i+1}(X)$.
As in~\cite[Corollary~4.14 and Theorem~4.15]{Bes99} we can prove the
following basic property of Coleman functions:
\begin{theorem}
  The sequence $0\to K \to \acol(X)\xrightarrow{\dd} \ocol^1(X)
  \xrightarrow{\dd} \ocol^2(X)$ is
  exact.
\end{theorem}
The proof is more or less the same as in the rigid case. First one
needs the algebraic version of Lemma~2.1 in~\cite{Bes99}.
\begin{lemma}
  Let $(E,\nabla)$ be a connection on $X$ with $E$ torsion free and let
  $f:E\to B$ be a map of coherent $\O_X$-modules with $B$ torsion free as
  well. Let $A\subset E$ be the kernel of $f$. Construct by induction
  a sequence of subsheaves:
  \begin{align*}
    A_0&=A\;,\\ A_{n+1}&=\Ker \left[A_n \xrightarrow{\nabla}
    E\otimes_{\O_X} \Omega_X^1 \to  (E\otimes_{\O_X} \Omega_X^1) /
    (A_n\otimes_{\O_X} \Omega_X^1)\right] \;.
  \end{align*}
  Let $A_\infty= \cap_{n=0}^\infty A_n$. Then 
  \begin{enumerate}
  \item $\nabla$ maps $A_\infty$ to $A_\infty \otimes_{\O_X}
    \Omega_X^1$.
  \item  The pair $\iny(A,\nabla):= (A_\infty,\nabla)$ is a
  connection.
\item  Any map of connections to $(E,\nabla)$ whose image lies in $A$
  factors through $\iny(A,\nabla)$.
  \end{enumerate}
\end{lemma}
Since the annalitification functor is exact on coherent sheaves it
follows
immediately that the analogue of Lemma~2.2 in loc.\ cit.\ is true in a
mixed algebraic, rigid analytic situation: If $U\in X^{\textup{an}}$ is a
rigid open subspace, then
$\nabla^{-\infty}(A|_U)=(\nabla^{-\infty}A)|_U$. This holds in
particular when $U$ is a disc in $X$. The same result holds with
respect to the maximal integral subconnection $E^{\textup{int}}$ of
the connection $(E,\nabla)$.
\begin{lemma}
  If a Coleman function $f\in \acol(X,\shF)$ vanishes on a disc, then
  it is identically $0$.
\end{lemma}
The proof is identical to the rigid proof found in Proposition~4.12 of
loc.\ cit. As in Corollary 4.14 of loc.\ cit., it follows that the
kernel of $\dd$ consists only of constants. The exactness of the
differentials at the one forms follows as in the proof of Theorem 4.15
in loc.\ cit.

As in~\cite[Definition 4.7]{Bes99} there is, for a morphism $f:X\to
Y$, a pullback map
\begin{equation*}
f^\ast:\acol(Y,\shF) \to \acol(X,f^\ast \shF)\;,
\end{equation*}
which gives as particular cases a ring homomorphism $f^\ast:\acol(Y)
\to \acol(X)$
and maps $f^\ast : \ocol^i(Y) \to \ocol^i(X)$ compatible with
differentials.

\begin{proposition}
  The association $U\mapsto \acol(U, \shF|_U)$ is a Zariski sheaf on $X$.
\end{proposition}
\begin{proof}
  The same proof as in~\cite[Proposition~4.21]{Bes99} works.
\end{proof}

In~\cite[Section~6]{Bes99} we defined an operator on
a subspace of
Coleman functions which we termed \emph{the $p$-adic
  $\db$-operator}. Here we recall this theory in the context of
Vologodski's theory. Everything works essentially without any
change. We define the subspace $\acoln_n(X,\shF)$ of $\acol(X,\shF)$
to consist of all Coleman functions that have a representative
$[M,\nabla,s,y]$ where $(M,\nabla)$ is a successive extension of at
most $n+1$ trivial connections (where here trivial means a direct sum
of any number of copies of $\O_X$ with its trivial connection). As in
loc.\ cit.\ this space is locally described by iterated integrals
involving at most $n$ iterated integrals. For an open $U\subset X$ we define
\begin{equation}
  \label{eq:ht}
  \htt_{\shF}(U):= \hdr^1(U/K)\otimes_K \shF(U)\;.
\end{equation}
This is a presheaf on $X$ and there are interesting obstructions for
gluing (see \cite[Corollary 6.9]{Bes99}).
The map $\db: \acoln_1(X,\shF) \to \htt_{\shF}(X)$ is defined as follows:
Suppose we are given a Coleman
function $[E,\nabla,s,y]$ where the connection $(E,\nabla)$ sits in a
short exact sequence
\begin{equation*}
  0\to E_1 \to E \to E_2\to 0\;,
\end{equation*}
such that $E_1$ and $E_2$ are trivial.
The projection of the $y_x$ give a compatible system of horizontal
sections of $E_2$. Since $E_2$ is trivial
this system comes from a
global horizontal section $y_2$ of $E_2$.
The connection $(E,\nabla)$ gives an extension class $[E]\in
\Ext_\nabla^1(E_2,E_1)$. The horizontal section $y_2$ is an element of
$\Hom_\nabla(\O_X,E_2)$. We can pullback the extension $[E]$ via $y_2$
to obtain $[E]\circ y_2\in
\Ext_\nabla^1(\O_X,E_1)$. The homomorphism $s$ restricts to $s_1 \in
\Hom(E_1,\shF)$ Since $E_1$ is trivial the
natural map
\begin{equation*}
  \Hom_\nabla(E_1,\O_X)\otimes \Hom(\O_X,\shF)\to \Hom(E_1,\shF)
\end{equation*}
is an isomorphism. Thus we may view $s_1$ as an element of the
left hand side. There is a product
$\Hom_\nabla(E_1,\O_X) \otimes \Ext_\nabla^1(\O_X,E_1) \to
\Ext_\nabla^1(\O_X,\O_X)$. Taking this product in the first coordinate
of $s_1$ with $[E]\circ y_2$ we obtain
\begin{equation*}
  \db (E,\nabla,s,y) :=
  -([E]\circ y)\circ s^\prime \in \Ext_\nabla^1(\O_X,\O_X) \otimes
  \Hom(\O_X,\shF)=\hdr^1(X/K) \otimes \shF(X)
\end{equation*}
(note that we changed the sign from loc.\ cit.\ and that
Proposition~6.8 there is actually true with the new sign).
As in loc.\ cit.\ Proposition~6.4 the $\db$ operator is independent of
the choice of the
short exact sequence in which $E$ sits and on the choice of a
representative abstract Coleman function. Also,
as in loc.\ cit.\ Proposition~6.5 this operator can be
described when $X$ is affine as sending 
$(\int \omega) \times s$ to $[\omega]\cdot s$, where $[\omega]$ is
the cohomology class of $\omega$.
\begin{proposition}\label{volddbsurjects}
  There is a short exact sequence,
  \begin{equation*}
    0\to \shF(X)\to \acoln_1(X,\shF)\xrightarrow{\db} \htt_{\shF}(X)\;,
  \end{equation*}
  which is exact on the right if $X$ is affine.
\end{proposition}
\begin{proof}
The proof is the same as in~\cite[Proposition~6.6]{Bes99}.
\end{proof}
A final observation regarding $\db$, which is obvious from the above
description, is the following.
\begin{lemma}\label{dercurve}
  If $l:\shF \to \shG$ is an $\O_X$-linear map and $F\in
  \acoln_1(X,\shF)$, then $l\circ F\in \acoln_1(X,\shG)$ and
  $\db (l\circ F)= (\id\otimes l)(\db F)$.
\end{lemma}
 From now on we will denote $\htt_{{\Omega_{X/K}^1}}$ simply by $\Ht$, as
this is the case that will be used almost exclusively in this work.

Suppose now that $X$ is a smooth variety over $\qpc$. Then, working
with models over finite extensions of $\Q_p$ and using functoriality,
it is clear that we can extend all of the above constructions to the
case $K=\qpc$. We can easily recover the field of definition of a
Coleman function as follows.
\begin{proposition}\label{2.9}
  Let $X/K$ be a smooth variety, $\shF$ a locally free sheaf defined
  over $K$, and let $\overline{X}$ and $\overline{\shF}$ be the extension of
  these objects to the algebraic closure of $K$. We assume that we
  have chosen a branch of the logarithm defined over $K$.  Let $f\in
  \acol(\overline{X},\overline{\shF})$ be a Coleman function and let
  $\sigma$ be an
  automorphism of  $\overline{K}$ over $K$. Then the function
  $f^\sigma$ defined by $f^\sigma(x)=\sigma(f(\sigma^{-1}(x)))$ is
  also a Coleman function on $\overline{X}$. If for any such $\sigma$
  we have $f^\sigma=f$, then in fact $f$ is an extension of scalars
  from a function in $\acol(X,\shF)$. 
\end{proposition}
\begin{proof}
It is obvious that $f^\sigma$ is a Coleman function. The second
statement is proved using minimal models for Coleman functions in the
same way that the sheaf property is proved.
\end{proof}

If $\pi:X\to Y$ is a finite covering of varieties over $\qpc$ and $f$ is a Coleman
function on $X$
we can construct the trace of $f$ along $\pi$ down to
$Y$. Surprisingly perhaps, this may not always be a Coleman function
on $Y$. We will now describe a very simple situation where one can
show that the trace is indeed a Coleman function. The problem seems
interesting and deserves further study, but here we limit ourselves to
a situation that suffices for our uses in this paper.
\begin{lemma}\label{tracecurve}
  Suppose $\pi:X\to Y$ is a finite covering and $\eta$ is a locally
  analytic one-form on $Y$ such that
  \begin{enumerate}
  \item  $\omega:= \pi^\ast \eta$ belongs to $\Ocoln_1(X)$, and
  \item There exists $\alpha\in \Ht(Y)$ such that $\pi^\ast
    (\alpha)=\db \omega$.
  \end{enumerate}
  Then in fact $\eta\in \Ocoln_1(Y)$ and $\db \eta = \alpha$.
\end{lemma}
\begin{proof}
The problem is (Zariski) local on $Y$, which we may therefore assume
affine. By Proposition~\ref{volddbsurjects} we can find $\eta' \in 
\Ocoln_1(Y)$ such that $\db \eta' = \alpha$. It follows that $\db
(\omega-\pi^\ast \eta')=0$ and therefore, again by
Proposition~\ref{volddbsurjects}, $\omega-\pi^\ast
\eta'= \pi^\ast (\eta-\eta')\in \Omega^1(X)$. It follows that
\begin{equation*}
  \eta-\eta'=\frac{1}{\deg \pi} \tr_\pi \pi^\ast (\eta-\eta') \in
  \Omega^1(Y)\;,
\end{equation*}
which implies the result.
\end{proof}
\begin{proposition}
  Suppose we have a diagram of finite maps $X'\xrightarrow{\pi'} X
  \xrightarrow{\pi}Y$ where the composition $X'\xrightarrow{\pi''} Y$
  is a Galois covering with Galois group $G$. Let $F\in \acol(X)$ with
  $\partial F\in \Ocoln_1(X)$ (Here we use $\partial$ instead of $\dd$
  so that the composed operator is $\dbd$, which is similar to the
  complex notation). Suppose there exists $\alpha\in \Ht(Y)$ such that
  \begin{equation*}
    (\pi'')^\ast \alpha 
    =\sum_{\sigma\in G} \sigma^\ast((\pi')^\ast \dbd F)\;.    
  \end{equation*}
  Then $\tr_\pi(F)\in \acol(Y)$ and $\dbd \tr_\pi(F)=\alpha/\deg(\pi')$.
\end{proposition}
\begin{proof}
We have $\tr_\pi(F)=\tr_{\pi''}((\pi')^\ast F)/\deg(\pi')$, so it
suffices to consider the case $\pi'=\id$, $\pi=\pi''$. Applying
Lemma~\ref{tracecurve} to $\eta=\dd \tr_\pi(f)$ gives the result.
\end{proof}

\section{The double index}
\label{sec:double}

In this section we recall the theory of the double index
from~\cite[Section~4]{Bes98b}. We must do a few things anew for the
algebraic theory used in this paper.

Let $K$ be a field of characteristic $0$. We
consider the field of Laurent series in the variable $z$ over $F$,
$\mer=K((z))$, the polynomial algebra over $\mer$ in the formal
variable $\log(z)$, 
$\alog:=\mer[\log(z)]$, and 
the module of differentials $\Olog:=\alog \cdot \dd z$. There is a formal
derivative  $\dd:\alog \to \Olog$ such that $\dd \log(z)=\dd z/z$ and it is
an easy exercise
in integration by parts to see that every form in $\Olog$ has
an integral in $\alog$ in a unique way up to a constant. We
distinguish in $\alog$ the subspace $\aloga=\mer + K\cdot \log(z)$
consisting exactly of all functions whose differential is in $\mer
\cdot \dd z$. To $F\in \aloga$ we can associated the residue of its
differential $\res \dd F\in K$. If $F\in \aloga$, then $F\in \mer $ if
and only if $\res \dd F=0$.
\begin{definition}
  The \emph{double index} is the unique anti-symmetric bilinear form
  $\pair{~,~}: \aloga \times \aloga \to K$ with the property that
  $\pair{F,G} = \res F \dd G$ whenever the left hand side has a meaning.
\end{definition}
The existence of the double index relies on a trivial linear algebra
lemma (Lemma~4.4 of \cite{Bes98b}). In loc.\ cit.\ it is computed in a
rigid analytic context and denoted $ind(F,G)$. Assume from now on that
$K$ is a complete subfield of $\C_p$ and that a branch of the log on
$K$ has been
chosen. The following lemma is
the algebraic analogue of Lemma~4.6 in loc.\ cit.
\begin{lemma}\label{pullindex}
  With a subscript $z$ or $w$ denoting the variable, let $\alpha:\mer_z
  \to \mer_w$ be given by $\alpha(z)=\sum_{k=n}^\infty a_k w^k$, with
  $a_n\ne 0$, and extend $\alpha$ in the obvious way to $\alog$ and
  $\Olog$. Then $\pair{\alpha(F),\alpha(G)}_w = n
  \pair{\alpha(F),\alpha(G)}_z$.
\end{lemma}

Suppose now that $K$ is either $\C_p$ or $\qpc$. For both of these
fields one has a ``Coleman integration theory'' of holomorphic
forms. This means that for any smooth connected variety $X/K$ there is
an integration map, $\omega\mapsto \int \omega$, from $\Omega^1(X)^{\dd=0}$ to
locally analytic $K$-valued functions on $X$ modulo the constant
functions, which is an inverse to the differential $\dd$ and which is
functorial with respect to arbitrary morphisms in the sense that for
such a morphism $f$ we have $\int (f^\ast \omega)= f^\ast  \int
\omega$. We further require that on $\PP^1$ we have $\int \dlog(z)=
\log(z)$ (in fact, it is not hard to show that the theory determines
uniquely a branch of the log for which this relation holds). Such a
theory over $K=\qpc$ exists as part of the more general theory
described in Section~\ref{sec:review}. For $\C_p$, or more generally
for complete subfields of $\C_p$ it follows from~\cite[Th\'eor\`eme
0.1]{Colm96}.

Let $X$ be a proper smooth curve over $K$ and suppose $F$ and $G$ belong to
$\acol(U)$, with $U\in X$ open, and such that $\dd F,\dd G\in
\Omega^1(U)$. Then we can compute the double index at every point
$x\in X$ in an obvious way, as the last lemma shows that the
computation will be independent of the choice of a local variable. The
index is non-zero only for $x\in X-U$. Since the sum of the
residues of $\dd F$ and $\dd G$ on $x\in X$ is $0$, the sum of all these
double indices does not change if we change $F$ or $G$ by a constant
and therefore we obtain a global pairing
\begin{equation}
  \label{eq:globalpair}
  \pair{\dd F,\dd G}_{\gl}:= \sum_{x\in X} \pair{F,G}_x\;.
\end{equation}
\begin{lemma}\label{seckindcase}
If both $\dd F$ and $\dd G$ are forms of the second kind, then
$\pair{\dd F,\dd G}_{\gl}=[\dd F]\cup [\dd G]$, where $[\dd F]$ and $[\dd G]$ denote the
cohomology classes in $\hdr^1(X)$ of the corresponding forms.
\end{lemma}
\begin{proof}
This follows from the usual formula for cup products on curves.
\end{proof}
In~\cite[Proposition~4.10]{Bes98b} we proved, under the assumption
that $X$ had good
reduction, an extension of this cohomological formula for the pairing,
where the forms are arbitrary, and in fact are even allowed to be defined
only on a certain rigid open
subset of $X$. For algebraic differentials of the third kind a similar
result was proved in the good reduction case by
Coleman~\cite[Theorem~5.2]{Col89} and in general by
Colmez~\cite[Th\'eor\`eme~II.4.2]{Colm96}. Here we will prove in
Theorem~\ref{gencolm} the analogue
of our result for meromorphic differentials but without assuming good
reduction, thus generalizing the result of Colmez.

\begin{proposition}\label{3.7}
  for any rational function $f\in K(X)$ and any meromorphic form
  $\omega$ we have $\pair{\omega,\dlog f}_\gl=0$.
\end{proposition}

We first need a few auxiliary results
\begin{lemma}\label{baseindex}
  On $\PP^1$ we have for any $a\in K$ that
  $\pair{\dlog(t),\dlog(t-a)}_{\gl}=0$.
\end{lemma}
\begin{proof}
  The non-trivial local indices are at $0$, $\infty$ and $a$. At
  $\infty$ we have
  \begin{equation*}
    \pair{\log(t),\log(t-a)}_\infty = \pair{\log(t),\log(t)}_\infty +
    \pair{\log(t),\log(1-\frac{a}{t})}_\infty = \log(1-\frac{a}{\infty})=0\;
  \end{equation*}
  while $ \pair{\log(t),\log(t-a)}_0=-\log(0-a)$
  and $ \pair{\log(t),\log(t-a)}_a=\log(a)$, so the result is clear.
\end{proof}
\begin{lemma}
  if $f:X\to Y$ is a finite map of curves, then $\pair{f^\ast
    \dd F,f^\ast \dd G}_\gl= \deg(f) \pair{\dd F,\dd G}_\gl$.
\end{lemma}
\begin{proof}
This follows immediately from the fact that if $f(x)=y$ with
multiplicity $n$, then we have $\pair{f^\ast
  F,f^\ast G}_x= n\pair{F,G}_y$ by Lemma~\ref{pullindex}.
\end{proof}
\begin{lemma}
  If $f:X\to Y$ is a finite map of curves and $\omega\in \Omega^1(K(Y))$,
  $\eta\in \Omega^1(K(X))$, then
  $\pair{f^\ast \omega,\eta}_{\gl} =\pair{\omega, \tr_f \eta}_\gl$.
\end{lemma}
\begin{proof}
We may assume that $f$ is a Galois covering with Galois group
$G$. Then we have
\begin{align*}
  \pair{f^\ast \omega,\eta}_\gl &= \frac{1}{|G|} \sum_{\sigma\in G}
  \pair{\sigma^\ast f^\ast \omega,\sigma^\ast \eta}_\gl =
  \frac{1}{|G|} \pair{ f^\ast \omega, \sum_{\sigma\in G}\sigma^\ast
    \eta}_\gl \\ &=
  \frac{1}{|G|} \pair{ f^\ast \omega, f^\ast \tr_f \eta}_\gl=
  \pair{\omega,\tr_f \eta}_\gl\;.
\end{align*}
\end{proof}

\begin{proof}[Proof of Proposition~\ref{3.7}]
The last lemma implies that $\pair{\omega,\dlog f}_\gl=\pair{\tr_f
\omega, \dlog t}_\gl$ on $\PP^1$. Now, $\tr_f \omega$ is a sum of a
form of the second kind on $\PP^1$ and a linear combination of forms
$\dlog(t-a_i)$. By Lemma~\ref{baseindex} we may assume that $\tr_f
\omega$ is of the second kind, and must therefore equal $\dd g$ for
some $g\in K(\PP^1)$. But $\pair{\dd g,\dlog t}_\gl= \sum_{x\in
\PP^1}\res_x (g\dlog t)=0$.
\end{proof}

\begin{definition}
  We define a map $\Psi': \Omega^1(K(X)) \to \hdr^1(X)$ by the
  condition that for any form of the second kind $\omega$ on $X$ we
  have
  \begin{equation*}
    \Psi'(\eta) \cup [\omega] = \pair{\eta,\omega}_\gl\;.
  \end{equation*}
\end{definition}
It follows immediately from Lemma~\ref{seckindcase} and
Proposition~\ref{3.7} that $\Psi'$ is well
defined and vanishes on log-differentials (i.e., of the form $\dlog f$). From
Lemma~\ref{seckindcase} we have $\Psi'(\eta)=[\eta]$ if $\eta$ is of
the second kind. The main result of this section is the following
theorem.
\begin{theorem}\label{gencolm}
  For any two meromorphic differentials $\omega$, $\eta\in
  \Omega^1(K(X))$ we have
  \begin{equation*}
    \pair{\eta, \omega}_\gl = \Psi'(\eta) \cup \Psi'(\omega)\;.
  \end{equation*}
\end{theorem}

To explain why this is a generalization of the results mentioned
above, we state the relation between the map $\Psi'$ and the logarithm
on the universal vectorial extension of the Jacobian of $X$.

Let $\T$ be the group of differentials of the third kind on $X$
and let $\T_l$ be the subspace of log differentials.

Recall from~\cite[Proposition~2.5]{Col-Gro89} that there exists a
commutative diagram with exact rows
\begin{equation*}
  \xymatrix{
    0 \ar[r] &\Omega^1(X)\ar@2{-}[d]  \ar[r]& \T/\T_l
    \ar[r]\ar[d]^{\Psi}& J \ar[r]\ar[d]^{\log_J} &0\\ 
    0\ar[r] &\Omega^1(X) \ar[r] &\hdr^1(X) \ar[r] &
    \hdr^1(X)/\Omega^1(X) \ar[r] & 0
  }
\end{equation*}
where $J$ is the Jacobian of $X$, the map $\T\to J$ sends a
differential of the third kind to its residue divisor and $\log_J$ is the
logarithm from $J$ to its Lie
algebra, which is isomorphic to $\hdr^1(X)/\Omega^1(X)$. The group
$\T/\T_l$ is the group of $K$ points of the universal vectorial
extension $G_X$ of $J$ and $\Psi$ is simply the logarithm for this
group. To be precise, we are working here with $K$ points, which is
fine for $K=\C_p$. For $K=\qpc$ we are taking the limit of the
corresponding map over all finite extensions of $\Q_p$. Note that for
a commutative algebraic group $H$ the logarithm
is characterized
as the unique locally analytic group homomorphism whose differential
at the identity element is the identity.

Since the map $\Psi'$ vanishes on $\T_l$ it induces a map $\Psi': \T/\T_l
\to \hdr^1(X)$.
\begin{proposition}\label{psiispsiprime}
  We have $\Psi'=\Psi$ on $\T/\T_l$.
\end{proposition}

To prove both the last proposition and our theorem we need to recall
the notion of differentiation of differential
forms with respect to a vector field (see~\cite[2.24]{War83}
or~\cite[Proposition 4.6]{Yek95}).
\begin{definition}
  Let
  $Y$ be a space ($Y$ could be a scheme or an analytic space over some
  field or a $C^\infty$ space) and let $\partial/\partial t$ be a vector
  field on $Y$. Then the operator of differentiation with respect to
  $\partial/\partial t$ is defined on $\Omega^i(Y)$ by
  \begin{equation}\label{yek}
    \frac{\partial}{\partial t} \omega = \dd
    (\omega|_{\frac{\partial}{\partial t}})+(-1)^i
    (\dd\omega)|_{\frac{\partial}{\partial t}}
  \end{equation}
  where $|_{\partial/\partial t}$ is the contraction operator.
\end{definition}
This operator commutes with the exterior differential $\dd$.
\begin{lemma}\label{derdoub}
  Let $\omega$ be a meromorphic form on $X$. Let $T$ be a
  smooth parameter variety and let $(\eta_t)_{t\in T}$ be a
  family of forms of the third kind on $X$ parametrized by $T$ (see
  Appendix~\ref{sec:univ} for the this notion) and
  $\partial/\partial t$ a vector field on
  $T$. Then $\frac{\partial}{\partial t} \eta_t $ is a form of the
  second kind, $\pair{\omega,\eta_t}_\gl$ is locally analytic in $t$ and
  \begin{equation*}
    \frac{\partial}{\partial t} \pair{\omega,\eta_t}_\gl=
    \pair{\omega, \frac{\partial}{\partial t} \eta_t}_\gl\;.
  \end{equation*}
\end{lemma}
\begin{proof}
The fact that the derivative is of the second kind is prove in the
Appendix.
We check the formula  at some $t_0\in T$.
We can find discs $D_i\subset X$, $i=1,\ldots, n$ and a ball $B$ around
$t_0$ such that for each $t\in B$ all the singular points of any of
the forms $\eta_t$ are contained in some $D_i$. By
Proposition~\ref{3.7} we may change $\omega$ by a linear combination
of log differentials, so we may assume $\omega$ has no residues inside
the $D_i$'s. By further reducing $D_i$ we may assume
that $\int \omega$ is a meromorphic function on each $D_i$. In such a situation,
the argument of the proof of Proposition~5.5 in~\cite{Bes98b} implies
that for each $t\in B$ we may
replace $\pair{\omega,\eta_t}_\gl$ with $\sum_i
\pair{\int \omega, \int \eta_t}_{e_i}$, where $e_i$ is an annulus
around $D_i$ 
and the local index around an annulus is the one defined
in~\cite[Proposition~4.5]{Bes98b}. In this case, each double index
$\pair{\int \omega,\int \eta_t}_{e_i}$ equals $\res_{e_i} \left((\int \omega)
\eta_t\right)$. By writing a Laurent expansion it is clear that each of these
expressions commute with differentiation.
\end{proof}
\begin{proof}[Proof of Theorem~\ref{gencolm}]
Let $\eta'$ be a form of the second kind representing
$\Psi'(\eta)$. Then, by the definition of $\Psi'$ we have
\begin{equation*}
  \Psi'(\eta) \cup \Psi'(\omega) = \pair{\eta', \omega}_\gl
\end{equation*}
so we must prove that 
\begin{equation}\label{eq:reqform}
\pair{\eta-\eta', \omega}_\gl=0
\end{equation}
 for any
$\omega$. Again by the definition of $\Psi'$ this is true for any
$\omega$ of the second kind. Since $\pair{\eta-\eta', \omega}_\gl$ is
linear in 
$\omega$ it suffices to prove that it vanishes on forms of the third
kind, and we already know that it vanishes on log-differentials by
Proposition~\ref{3.7}. We thus get an additive map
\begin{equation*}
  \T/\T_l \to K,\quad \omega\mapsto  \pair{\eta-\eta', \omega}_\gl\;,
\end{equation*}
and we must show that it is the $0$ map. Since $\T/\T_l=G_X(K)$, where
$G_X$ is the universal vectorial extension of $J$, it suffices to show
that the derivative of this map at the identity is $0$. But this is
clear by Lemma~\ref{derdoub} because the derivative of a family of
forms of the third kind is a form of the second kind $\omega$, for
which \eqref{eq:reqform} holds.
\end{proof}
\begin{proof}[Proof of Proposition~\ref{psiispsiprime}]
Since $\Psi'$ is clearly additive, it suffices to show that it is
analytic near the origin and that its differential at $0$ is the
identity map when $\Lie(G_X)$ is identified with $\hdr^1(X)$. This is
equivalent to showing that for any form $\omega$ of the second kind
the map $G_X\to K$, given by
\begin{equation}\label{showanal}
  \eta\mapsto  \Psi'(\eta) \cup [\omega] = \pair{\eta,\omega}_\gl\;,
\end{equation}
is analytic near the origin and its differential at $0$ is given by
cupping with $[\omega]$. By Lemma~\ref{loclift} we can find an
algebraic family of forms $\eta_t$ for $t$ in a neighborhood of $0\in G_X$
such that $\eta_t$ represents $t$. This already shows that
\eqref{showanal} is analytic near $0$. The formula for the
differential at $0$ follows easily from Lemma~\ref{derdoub} and
Corollary~\ref{A5}.
\end{proof}
The following corollary is the only result needed in the rest of the text.
\begin{corollary}\label{doubleisproj}
   Suppose $\omega$ is a form of the second kind while $\eta\in \T$. Then
   $\pair{\omega,\eta}_{\gl}= [\omega]\cup \Psi(\eta)$.
\end{corollary}

The following result is needed in~\cite{Bes99a}. We give it here since
it uses the techniques of this section.
\begin{proposition}\label{ColemanisColmez}
  Let $X$ be a curve over $\C_p$ with good reduction. Then the Coleman
  and Colmez integrals of algebraic forms on $X$ coincide.
\end{proposition}
\begin{proof}
For forms of the second kind this follows from~\cite{Col85} as it is
shown there that the Coleman integral can be obtained by pullback from
an abelian variety, which is how the Colmez integral is defined. By
linearity of both integrals it
remains to consider a form of the third kind $\eta$. Let $F_1$ and
$F_2$ be the Colmez (respectively Coleman) integral for $\eta$. Since
both $F_1$ and $F_2$
differentiate to $\eta$, $F_1-F_2$ is a locally
constant function on $X$. For a divisor $D$ of 
degree $0$ on $X$ let $\beta(D)$ be $F_1-F_2$ evaluated at
$D$. The function $\beta$ is locally constant on the points in
$D$. It is furthermore additive. If we show that it vanishes on
principal divisors it will give a locally constant additive function
on the Jacobian and will have to vanish. So suppose $D=(f)$ for a
rational function $f$. Both Coleman and Colmez integrals of $\dlog(f)$
are $\log(f)$. It follows from
Proposition~\ref{3.7} that $\sum_{x\in X} \pair{\log(f),F_1}_x=0$ and
the same follows for $F_2$ by~\cite[Corollary~4.11]{Bes98b} combined
with the analysis in the proof of Proposition~5.5 in~\cite{Bes98b}
that shows we may replace local indices on annuli by local indices at
points. But $\sum_{x\in X} \pair{\log(f),F_1-F_2}_x= \beta(D)$ so the
proof is complete.
\end{proof}

\section{Metrized line bundles over the $p$-adics}
\label{sec:metrized}

In this section it is most convenient to work over an algebraically
closed field, so we will always work over $\qpc$: All varieties
will be smooth and connected over $\qpc$, and differential
forms and de Rham cohomology will also be taken over $\qpc$. We also
fix a branch of the logarithm $\log:\qpc \to \qpc$.

In classical Arakelov theory one of the basic notions is that of a
metrized line bundle: A line bundle endowed with a metric. Often one
does not use the absolute value of the section but rather the log of
this absolute value. Experience has shown that the replacement for the
log of an absolute value in the complex case is simply the p-adic
logarithm in the p-adic case. This suggests the following definition.
\begin{definition}
  Let $X$ be a smooth variety and
  let $\L$ be a line bundle over $X$. A \emph{log function} on $\L$ is a
  function $\log_\L\in \acol(\Lt)$, where $\Lt :=\tot(\L)-\{0\}$,
  $\tot(\L)$ being the total space of $\L$,
  such that the following conditions are satisfied:
  \begin{enumerate}
  \item For any $x\in X$, any $v$ in the fiber $\L_x$ and $\alpha\in
    \qpc$ one has 
    $\log_\L(\alpha v)= \log(\alpha)+\log_\L(v)$.
  \item $\dd \log_\L\in \Ocoln_1(\Lt)$
  \end{enumerate}
  A function satisfying only the first condition will be called a
  \emph{pseudo log function}. 
  A line bundle together with a log function will be called a \emph{metrized
  line bundle} (the more natural name of a log line bundle is rejected
  because of possible confusion with the theory of log schemes). There
  is an obvious notion of isometry of metrized line bundles.
\end{definition}
Clearly the first condition implies that if $s$ is a section of $\L$
and $f$ is a rational function on $X$, then $\log_\L(fs)=\log(f) +
\log_\L(s)$. Adding a constant to a log function one obtains a new log
function. This operation will be called scaling.
\begin{remark}
The second condition on a log function seems a bit arbitrary. Clearly,
some condition should exist to make the log function amenable to the
tools of $p$-adic analysis. A weaker condition could be that for one,
hence any, invertible section $s$ of $\L$ in a Zariski open $U$ the function
$\log_\L(s)$ on $U$ is a Coleman function. This is clearly implied by
our stronger condition. One reason for this condition is that we can
use the theory of the $p$-adic $\db$ operator to define a kind of
curvature for our log function. Another reason is that the canonical
Green function we will define can be characterized, as we will see in
Proposition~\ref{charofgr}, by different means, and satisfies this
property. As a
consequence, log functions admissible with respect to this Green
function, in a sense we will define, have this property and these are
the only log functions we are interested in.
\end{remark}
\begin{definition}
  If $\L$ and $\mathcal{M}$ are metrized line bundles then $\L\otimes
  \mathcal{M}$ 
  has a canonical log function $\log_{\L}\otimes \log_{\mathcal{M}}$
  defined by 
  \begin{equation*}
   (\log_{\L}\otimes \log_{\mathcal{M}}) (s\otimes t) = \log_\L(s)+
    \log_{\mathcal{M}}(t)\;.
  \end{equation*}
\end{definition}

Exactly as in the complex case, it turns out one can usually associate
to a metrized line bundles a kind of curvature. Let $\L$ be a metrized
line bundle on $X$. Since
$\dlog_\L\in \Ocoln_1(\Lt)$ by assumption we can take its $\db$ operator. We
use the notation $\dbd \log_\L\in \Ht(\Lt)$ for the result instead of using
$\db \dd$ to make the notation similar to the complex one.
\begin{proposition}\label{excurve}
  Suppose $X$ is proper and let $\pi: \Lt \to X$ be the projection.
  \begin{enumerate}
  \item Suppose that $ch_1(\L)\in \im\left(\cup: \Ht(X)\to
    \hdr^2(X)\right)$. Then there exists a unique $\Curve(\L)\in \Ht(X)$
    such that $\pi^\ast \Curve(\L)=\dbd \log_\L$. The element
    $\Curve(\L)$ is called the \emph{curvature form} of the metrized
    line bundle $\L$ and it satisfies the relation $\cup
    \Curve(\L)=ch_1(\L)$.
  \item Conversely, suppose that $\alpha\in \Ht(X)$ satisfies $\cup
    \alpha = ch_1(\L)$. Then there exist a log function $\log_{\L}$ on
    $\L$ such that $\Curve(\L)=\alpha$.
  \end{enumerate}
\end{proposition}
\begin{proof}
Consider first the case of the trivial line bundle $\L$ over $X$,
where we have
$\Lt=X\times (A^1-\{0\})$. The log function is then given by
$\log_{\L}(x,t)= \log_{\L}(1)(x)+\log(t)$. Since $\dbd(\log(t))=0$ it
follows that in this case $\dbd(\log_{\L})=\pi^\ast (\dbd
\llog(1))$. It is also clear that $\dbd \llog(1))$ is characterized by
this equation. Suppose now that $\L$ is arbitrary and that
$\mathcal{U}=\{U_i\}$ is a covering of
$X$ over which $\L$ is trivialized, so that we have nonvanishing
sections $s_i$ over $U_i$. We obtain a $0$-\cech\ cocycle $i\mapsto
\alpha_i:=\dbd(\llog(s_i))$, since
$\llog(s_i)-\llog(s_j)=\log(s_i/s_j)$ on $U_i\cap U_j$ and its $\dbd$ is
$0$. If this cocycle comes from $\alpha\in \Ht(X)$ then $\pi^\ast
\alpha$ equals $\dbd \llog$ on $\pi^{-1}(U_i)$ for each $i$ and thus
on $\Lt$, proving the existence of the curvature. To show that we can
indeed find such an $\alpha$ we now use the criterion
of~\cite[Corollary~6.9]{Bes98a}. According to it,
$\alpha$ exists if the image of the cocycle $i\mapsto \alpha_i$ under
the map $\ninv$ of loc.\ cit.\ Definition~6.7, is in the image of
$\cup$. It thus suffices to prove that this image is $ch_1(\L)$. The
map $\ninv$ is computed as follows: For each $i$ we need to find
$\omega_i \in \Ocoln_1(U_i)$ such that $\db(\omega_i)=\alpha_i$ and
then consider the cocycle $ij\mapsto \omega_i-\omega_j$. But in our
case we can clearly take $\omega_i=\dllog(s_i)$ and
$\omega_i-\omega_j=\dlog(s_i/s_j)$. This last cocycle is well known to
represent $ch_1(\L)$ so the first assertion is proved, except for
uniqueness. This follows because clearly the map $\Ht(X)\to \Ht(U)$,
when $U$ is open in $X$, is injective.

Suppose now that we are given $\Omega=\sum \theta^i\otimes \omega^i\in
\Ht(X)$ such that $\cup \Omega=ch_1(\L)$. Notice that $\dd \omega_i=0$
since the $\omega_i\in \Omega^1(X)$ and $X$ is proper. We want to
construct a log
function on $\L$ whose curvature is $\Omega$. The first step is to
write explicit cocycles representing the $\theta^i$ and then $\cup
\Omega$. We do this using the affine covering
$\mathcal{U}=\{U_j\}$. With respect to this covering we can write
\begin{equation*}
  \theta^i = ((\eta_j^i\in \Omega^1(U_j)),(f_{jk}^i\in \O(U_{jk})))
  ,\; \eta_j^i-\eta_k^i=\dd f_{jk} \text{ on } U_{jk}\;,
\end{equation*}
Where $U_{jk}=U_j\cap U_k$.
We make the comparison $\cup \Omega = ch_1(\L)$ in $H^2(X,F^1
\Omega^\bullet)$ instead of in the full $\hdr^2(X/K)$. Since $X$ is
proper this forms a subspace and both sides of the equation belong to
this subspace. Using again \cech\ cocycles we can write
\begin{equation*}
  H^2(X,F^1 \Omega^\bullet)=
  \frac{ \{ (\chi_j\in \Omega^2(U_j), \zeta_{jk}\in
    \Omega^1(U_{ij}),\chi_j-\chi_k =\dd \zeta_{jk},
    \zeta_{jk}-\zeta_{jl}+\zeta_{kl}=0\}}{\{(\dd
    \alpha_j,\alpha_j-\alpha_k), \alpha_j \in \Omega^1(U_j)\}}.
\end{equation*}
With this representation we have
\begin{align*}
  \cup \Omega &= (\sum_i \eta_j^i \wedge \omega^i,\sum_i f_{jk}
  \omega^i)\\ \intertext{and}
  ch_1(\L) &= (0,\dlog g_{jk})\;,
\end{align*}
where $g_{jk}=s_j/s_k$.
Therefore, the condition $\cup \Omega = ch_1(\L)$ spells out as
\begin{equation}
  \label{eq:twoconds}
  \begin{split}
    \sum_i \eta_j^i \wedge \omega^i &= \dd \alpha_j\\
    \sum_i f_{jk} \omega^i &= \dlog g_{jk}+\alpha_j-\alpha_k\;,
  \end{split}  
\end{equation}
for some $\alpha_j\in \Omega^1(U_j)$.

To define the log function it suffices to define $\llog(s_j)$.
By assumption, 
\begin{equation*}
  \dbd \llog(s_j)=\Omega|_{U_j}= \sum_i [\eta_j^i]\otimes \omega^i\;.  
\end{equation*}
This implies that for some choices of $\gamma_j\in \Omega^1(U_j)$ and
Coleman integrals $H_j^i=\int \eta_j^i$ we should have
\begin{equation*}
  \dllog s_j=\sum_i H_j^i \omega^i+ \gamma_j\;.
\end{equation*}
We further need to have $\llog (s_j) - \llog(s_k) = \log(g_{jk})$ and
differentiating this we get
\begin{equation*}
  \sum_i (H_j^i-H_k^i) \omega^i + \gamma_j-\gamma_k = \dlog g_{jk}
\end{equation*}
(recall that $\dd\omega_i=0$).
Since $H^1(\mathcal{U},\qpc)=0$ we can arrange the constants of
integration in such a way that $H_j^i-H_k^i=f_{jk}^i$. We then get
\begin{equation*}
  \sum_i (f_{jk}^i \omega^i + \gamma_j-\gamma_k) = \dlog g_{jk}
\end{equation*}
and this can be arranged by taking $\gamma_j=-\alpha_j$ in view of
\eqref{eq:twoconds}. So set
\begin{equation*}
  \delta_j=\sum_i H_j^i \omega^i- \alpha_j
\end{equation*}
with the choice of $H_j^i$ discussed before. We notice that these
$\delta_j$ are closed forms. Indeed, since $\dd \omega^i=0$ we have
\begin{equation*}
  \dd \delta_j =\sum_i \eta_j^i \wedge \omega^i- \dd \alpha_j=0
\end{equation*}
by \eqref{eq:twoconds}. Define now a Coleman form on $\Lt$ as follows:
Choose an isomorphism $\pi^{-1} U_j\isom A^1\times U_j$ in such a way
that $s_j$ corresponds to the
section $1$ and define the form there by $\pi^\ast \delta_j+
\dlog(t)$. It now follows that these forms are closed and that they
glue to give a closed Coleman form $\delta\in \Ocoln_1(\Lt)$ which we
can then integrate to obtain our required $\llog$.
\end{proof}

The behavior of log functions with respect to pullbacks is given by
the following obvious result.
\begin{proposition}
  Let $\L$ be a line bundle on $Y$ with a log function $\llog$ and let
  $f:X\to Y$ be a morphism. Let $\L'=f^\ast \L$ and consider the map
  $\tilde{f}:{\L'}^\times \to \Lt$. Then $\tilde{f}^\ast \llog $ is a log
  function on $\L'$ whose curvature is $f^\ast \Curve(\L)$.
\end{proposition}
We want to consider the behavior of log functions with respect to
norms.
Suppose $\pi:X\to Y$ is a finite covering, $\L$ is a line bundle on $X$ with
a log function $\llog$. the norm of $\L$ to $Y$, $\norm_\pi \L$,
acquires a natural pseudo log function $\norm_\pi \llog$ as follows:
Its fiber over 
$y\in Y$ is $\otimes_{\pi(x)=y} \L_x^{\otimes{n_{\pi,x}}}$, with $n_{\pi,x}$
the multiplicity of $\pi$ at $x$. Since each of the fibers $\L_x$ has a
log function the tensor product has one as well.
\begin{lemma}
  If $\L$ has a pseudo log function $\llog$ and the induced pseudo log
  function on $\L^{\otimes n}$ is a log function, then so is $\llog$.
\end{lemma}
\begin{proof}
  This is clear once taking a \cech\ covering. If $s$ is an invertible
  section of $\L$ on $U$, then $s^{\otimes n}$ is an invertible
    section for $\L^{\otimes n}$, and if $\log_{\L^{\otimes n}}
    (s^{\otimes n})= n \llog(s)$ is a Coleman function of the right
    type, then so is $\llog(s)$.
\end{proof}

\begin{proposition}\label{5.5}
  Suppose we have a diagram of finite maps $X'\xrightarrow{\pi'} X
  \xrightarrow{\pi}Y$ where the composition $X'\xrightarrow{\pi''} Y$
  is a Galois covering with Galois group $G$. Let $\L$ be a line
  bundle on $X$ with log function $\llog$. Suppose there exists
  $\alpha\in \Ht(Y)$ such that
  \begin{equation*}
    (\pi'')^\ast \alpha 
    =\sum_{\sigma\in G} \sigma^\ast((\pi')^\ast \Curve(\llog))\;.    
  \end{equation*}
  Then $\norm_\pi(\llog)$ is a log function on $\norm_\pi \L$ and
  its curvature is $\alpha/\deg(\pi')$.
\end{proposition}
\begin{proof}
We have a natural isomorphism
\begin{equation*}
  \norm_{\pi''} (\pi^{\prime\ast} \L)\isom (\norm_\pi \L)^{\otimes\deg \pi'}
\end{equation*}
which is compatible with pseudo-log functions. By the previous lemma
the proposition reduces to the case where $\pi'=\id$ and $\pi=\pi''$ is
Galois. Let $\L'= \norm_\pi \L$ and let
\begin{equation*}
  \L'' :=\pi^\ast \norm_\pi \L \isom \bigotimes_{\sigma\in G} \sigma^\ast \L\;.
\end{equation*}
The last isomorphism is in fact an isometry. The log function on $\L''$
has curvature
\begin{equation*}
  \Curve(\log_{\L''})=
  \sum_{\sigma\in G}\sigma^\ast \Curve(\llog) = \pi^\ast \alpha\;.
\end{equation*}
Let  $\pi_Y:\L^{\prime \times} \to Y$, $\pi_X:\L^{\prime\prime\times} \to X$,
be the projections. The map $\tilde{\pi}: \L^{\prime\prime\times} \to
\L^{\prime\times}$ is finite. We have $\tilde{\pi}^\ast \dd \log_{\L'}=\dd
\log_{\L''}$ and
\begin{equation*}
  \db \dd \log_{\L''} = \pi_X^\ast \Curve(\log_{\L''}) = \pi_X^\ast
  \pi^\ast \alpha = \tilde{\pi}^\ast \pi_Y^\ast \alpha\;.
\end{equation*}
It follows from Lemma~\ref{tracecurve} that $\dd \log_{\L'}\in
\Ocoln_1(\L^{\prime\times})$ and $\db \dd \log_{\L'}=\pi_Y^\ast \alpha$ so the
curvature of $\log_{\L'}$ is $\alpha$.
\end{proof}

\section{The almost canonical Green function}
\label{sec:canonical}

We now construct the almost canonical Green function on a complete
non-singular curve $X$ of positive genus $g$ over $\qpc$. Taking the hint from
classical Arakelov
theory we define instead an almost canonical log function on
$\O(\Delta)$ on $X\times X$. Here almost canonical means canonical up
to scaling. The Green function is then simply $\log_{\O(\Delta)}(1)$
and it is therefore defined as a function of two variables up to a
constant.

We fix splitting $\hdr^1(X)=W\oplus \Omega^1(X)$. This type of
splitting occurs in the theory of $p$-adic height pairings.

Let $p:\O(\Delta)^\times \to X\times X$ and $\pi_1$, $\pi_2:X\times X\to
X$ be the obvious projections.
\begin{definition}\label{cancurve}
  We define elements $\mu\in \Ht(X)$ and $\Phi\in \Ht(X\times X)$ as
  follows: Fix a basis $\{\omega_1,\ldots,\omega_g\}$ of
  $\Omega^1(X)$. Let $\{\bo_1,\ldots, \bo_g\}\subset W$ be a dual basis
  with respect to the cup product
  (i.e., $\tr (\bo_i \cup \omega_j) =\delta_{ij}$). Then we set
  \begin{align*}
    \mu&=\frac{1}{g} \sum_{i=1}^g \bo_i \otimes \omega_i\in \Ht(X)\;,\\
    \Phi &= \pi_1^\ast \mu + \pi_2^\ast \mu
    - \sum_{i=1}^g \left(\pi_1^\ast \bo_i \otimes  \pi_2^\ast\omega_i
      +\pi_2^\ast \bo_i \otimes  \pi_1^\ast\omega_i
    \right) \in \Ht(X\times X)\;.
  \end{align*}
\end{definition}
\begin{lemma}
  We have $\cup \Phi= ch_1({\O(\Delta)})=cl(\Delta)$.
\end{lemma}
\begin{proof}
We have
\begin{equation*}
  \cup \Phi = \pi_1^\ast (\cup \mu) + \pi_2^\ast (\cup \mu)
    - \sum_{i=1}^g \left(\pi_1^\ast \bo_i \cup  \pi_2^\ast\omega_i
      +\pi_2^\ast \bo_i \cup  \pi_1^\ast\omega_i
    \right)\;.
\end{equation*}
To prove that this is the cohomology class of the diagonal it suffices
to show that $\tr((\cup \Phi)\cup \phi)= \tr \Delta^\ast \phi$ for any
$\phi\in \hdr^2(X\times X)$. By K\"unneth we can write such $\phi$
as a combination of forms of the following three forms: $\pi_i^\ast
(\cup \mu)$, $\pi_1^\ast \bo_i \cup \pi_2^\ast \omega_j$ or
$ \pi_2^\ast \bo_j \cup \pi_1^\ast \omega_i $. Checking the formula in
each of these three cases is straightforward.
\end{proof}

From the above lemma and Proposition~\ref{excurve} we obtain the
existence of a log function
on ${\O(\Delta)}$ with curvature $\Phi$. In the classical case the analogous
relation already suffices to characterize the metric up to a
constant. In our case however, this is not sufficient yet, since we
can always modify our log function by $p^\ast \int \phi$ for
$\phi\in \Omega^1(X\times X)$. Again by K\"unneth such a form $\phi$
can be written as $\sum_{i=1}^2 \pi_i^\ast \phi_i$ with $\phi_i\in
\Omega^1(X)$.

For any choice of a log function $\log_{\O(\Delta)}$ with curvature
$\Phi$ we can
define a corresponding Green function $G$ by $G=\log_{\O(\Delta)}(1)$,
where $1$
is the canonical section of $\O(\Delta)$. 
\begin{definition}\label{odcan}
  For any divisor $D=\sum n_i P_i$ on $X$ we define the Green function
  for $D$ as $G_D=\sum n_i G(P_i,\bullet)$. We define the canonical
  log function on $\O(D)$ by the condition $\log_{\O(D)}(1) = G_D$.
\end{definition}
We can
alternatively express this log function as follows: Suppose again $D=\sum
n_j P_j$. For $P\in X$ let
$i_P:X\to X\times X$ be the map $i_P(x)=(P,x)$. Since
$\O(P)=i_P^\ast(\O(\Delta))$ we have
\begin{equation*}
  \O(D)=  \otimes (i_{P_j}^\ast {\O(\Delta)})^{\otimes n_j}
\end{equation*}
and in this way $\O(D)$ inherits a log function from the log function
$\log_{\O(\Delta)}$.

Notice that the Green function for a divisor of degree $0$ is
determined by what was already done without fixing $\log_{\O(\Delta)}$ any
further. In fact, for principal divisors it is what one can expect.
\begin{proposition}\label{chidef1}
  We have $G_{(f)}=\log(f)+Const$.
\end{proposition}
\begin{proof}
This is equivalent to the following statement: The function $f$
determines an isomorphism $\O((f))\isom \O_X$ and this isomorphism is
an isometry up to a constant. To prove this, consider the map $f\times
\id: X\times X \to \PP^1\times X$. Choose any $x_0\in X$ and consider
the line bundle $\L={\O(\Delta)}\otimes \pi_1^\ast (\O(x_0))^{-1}$ with its induced
log function. We have
\begin{equation*}
  \Curve (\L)= \Curve({\O(\Delta)})- \pi_1^\ast \mu = \pi_2^\ast \mu
    - \sum_{i=1}^g \left(\pi_1^\ast \bo_i \otimes  \pi_2^\ast\omega_i
      +\pi_2^\ast \bo_i \otimes  \pi_1^\ast\omega_i
    \right)\;.
\end{equation*}
We want to show that this curvature has a trace via $f\times
\id$ and compute this trace. We have a diagram $X'\xrightarrow{\pi} X
\to \PP^1$ where the composed map $X'\to \PP^1$ is Galois, say with
Galois group $G$. This is then also true with respect to the base
change $X'\times X \to \PP^1 \times X$. Let $p_i$, $i=1,2$ be the
projections from $X'\times X$ to its factors. We now have
\begin{multline*}
  \sum_{\sigma\in G} (\sigma\times \id)^\ast (\pi\times \id)^\ast
  \Curve(\L)\\ =|G| p_2^\ast \mu-\sum_{i=1}^g\left(\sum_{\sigma\in G} p_1^\ast
  \sigma^\ast \pi^\ast \bo_i \otimes p_2^\ast \omega_i+\sum_{\sigma\in
  G}   p_1^\ast \bo_i  \otimes  p_2^\ast \sigma^\ast \pi^\ast \omega_i
  \right) = |G| p_2^\ast \mu\;.
\end{multline*}
The two sums in the brackets are $0$ because $\sum  \sigma^\ast
\pi^\ast \eta$, with $\eta=\omega_i$ or $\bo_i$, is a pullback from
$\PP^1$ and $\Omega^1(\PP^1)=\hdr^1(\PP^1)=0$.
Proposition~\ref{5.5} now implies that the 
induced quasi log function on $\L':= \norm_{f\times \id} \L$ is a
log function and its curvature is a multiple of $\tilde{\pi}_2^\ast
\mu$, where $\tilde{\pi}_2:\PP^1\times X \to X$ is the projection on the
second factor. Consider now $\L''=\L'\otimes (i_0^\ast \L')^{-1}$, where
$i_0:\PP^1\times X \to \PP^1\times X$ is given by $i_0(a,x):=(0,x)$.
We give $\L''$ the induced log function. The curvature of this
log function is $0$, which determines it up to the integral of an
element of $\Omega^1(\PP^1\times X)\xleftarrow{\sim \tilde{\pi}_2^\ast}
\Omega^1(X)$. Now, the restriction of $\L''$ to $0\times X$ is canonically
trivial with trivial log function. By what was said before the
conditions of trivial curvature and triviality on the restriction to
$0\times X$ determines the
log function uniquely. On the other hand, the restriction of $\L''$ to
$\infty \times X$ is up to scaling $\O((f))$ with its log
function. There is an isomorphism $\L'' \isom
\O_{\PP^1\times X}$. Indeed, suppose $f(x_0)=a_0$. The divisor corresponding to
$\L''$ is the graph of $f$ plus $a_0\times X$ minus the
pullback of the zero divisor of $f$. This is exactly minus the divisor of the
function 
\begin{equation*}
  F(a,x):=\left(\frac{1}{a} - \frac{1}{f(x)} \right) \cdot
  \frac{a}{a_0(a-a_0)}=
  \left(1 - \frac{a}{f(x)} \right) \cdot
  \frac{1}{a_0(a-a_0)}\;.
\end{equation*}
Thus, the trivial log function on $\O_{\PP^1\times X}$ induces a log
function on $\L''$ whose curvature is still $0$. Since $F(a,x)$
is constant for $a=0$ it is clear that this log function is trivial
after scaling. Thus, after scaling we obtained the same log function
as before on $\L''$. On the other hand, $F(\infty,x)^{-1}$ is a
constant multiple of $f(x)$, so up to scaling $f$ indeed induces an
isometry.
\end{proof}
\begin{corollary}
  An isomorphism $\O(D)\to \O(D\pr)$ is an isometry up to constant.
\end{corollary}
\begin{definition}\label{ilogdef}
  We define the homomorphism $\ilog:\qpc(X)^\times \to \qpc$ by
  $\ilog(f)=\log(f)-G_{(f)}$.
\end{definition}
The $\ilog$ character is the $p$-adic analogue of the integral of the
norm of a section. It will thus be used to associate the ``infinite''
fibers to the Arakelov divisor of a section of a line bundle.
\begin{definition}\label{ilogl}
  Let $\L$ be a line bundle on $X$. The log function $\log_\L$ is
  called admissible (with respect to the Green function $G$), if for each
  section $s$ of $\L$ the function $\log_\L(s) - G_{\Div(s)}$ is
  constant. This constant will be denoted $\ilog_\L(s)$.
\end{definition}
\begin{lemma}
  The condition for being admissible need only be checked on a single
  section. The function $\ilog_\L$ satisfies the relation $\ilog_\L(fs)
  = \ilog(f)+ \ilog_\L(s)$ when $f\in \qpc(X)^\times$.
\end{lemma}
\begin{proof}
  immediate from Proposition~\ref{chidef1}
\end{proof}
\begin{lemma}
  Any isomorphism between admissible metrized line bundles is an
  isometry up to scaling.
\end{lemma}
\begin{proof}
  If $T:\L_1 \to \mathcal{M}$ is an isomorphism let $s\in \L(C)$.
  Since $s$ and $T(s)$ have the same divisor admissibility implies
  that $\log_{\L}(s) - \log_{\mathcal{M}}(T(s))$ is a constant. Since
  any other section is obtained from $s$ via multiplication by a
  rational function it is easy to see that this constant is
  independent of $s$.
\end{proof}
For any divisor $D$ the log function we defined on the line bundle
$\O(D)$ is clearly admissible. Note also that the log function
$\log_\L \otimes \log_{\mathcal{M}}$ on $\L\otimes \mathcal{M}$ 
is admissible if both $\log_\L$ and $\log_{\mathcal{M}}$ are.

As a first step towards removing the degrees of freedom in the
definition of $\log_{\O(\Delta)}$ we can
assume that it is symmetric. Indeed, since there is a canonical
isomorphism between ${\O(\Delta)}$ and $\sigma^\ast {\O(\Delta)}$, where
$\sigma(x,y)=(y,x)$ we can consider, for any log function $\log_{\O(\Delta)}$ as
above the function $(\log_{\O(\Delta)}+\sigma^\ast \log_{\O(\Delta)})/2$. This still has the
same curvature form since $\Phi$ is symmetric. A symmetric log
function is determined up to a function of the form $\sum_{i=1}^2
p^\ast \pi_i^\ast \int \phi$ with $\phi\in \Omega^1(X)$.

It is well known that $\Delta^\ast \O(-\Delta)$ is canonically
isomorphic to $\bfo_X$ - the canonical bundle on $X$. Thus, given a
choice of a log function on ${\O(\Delta)}$ the line bundle $\bfo_X$ inherits a
canonical log function. It therefore make sense to impose the
condition that this log function is admissible with respect the Green
function the same log function defined. This allows us to determine a
canonical log function on ${\O(\Delta)}$, hence a canonical Green function, up
to a constant, as we will see in the next theorem.
\begin{theorem}\label{canG}
  There exist a unique up to constant symmetric log function on
  $\O(\Delta)$ with curvature $\Phi$ and such that the pulled back
  log function on $\bfo_X$ is admissible with respect to the
  induced Green function.
\end{theorem}
\begin{proof}
We interpret the admissibility condition as follows: for $P\in X$ let
$i_P:X\to X\times X$ be the map $i_P(x)=(P,x)$. We know that
$\bfo_X \isom \O(D)$ for some divisor $D=\sum n_j P_j$. Since
$\O(P)=i_P^\ast(\O(\Delta))$ we can write the isomorphism as
\begin{equation}\label{6.95}
  (\Delta^\ast {\O(\Delta)})^{-1} \isom \otimes (i_{P_j}^\ast {\O(\Delta)})^{\otimes n_j}\;.
\end{equation}
For each choice of a log function both sides inherit a log function
and we need to find one for which the isomorphism is an isometry up to
a constant. First we notice that both sides have the same
curvature. This is because $\Delta^\ast \Phi=(2-2g) \mu$ while
$i_P^\ast \Phi =\mu$ and $\deg(D)=2g-2$.
(this is analogous to the classical theory). Thus, the differential of
the difference of the log functions equal $\pi^\ast \omega$ for some
$\omega\in \Omega^1(X)$. On the other hand, if we modify the log
function by $\sum_{i=1}^2 p^\ast \pi_i^\ast \int \phi$ the
differential of the log function on the left hand side of \eqref{6.95}
changes by $-2 \phi$ while on the right hand side it changes by
$\deg D \cdot \phi = (2g-2) \phi$. Since $g>0$ by assumption we can solve
this equation 
uniquely to obtain $\dd\log_{\O(\Delta)}$ uniquely, hence
$\log_{\O(\Delta)}$ uniquely up to a constant.
\end{proof}

\begin{corollary}\label{dbdG}
  We have $\dbd G=\Phi|_{X\times X-\Delta}$.
\end{corollary}
In the course of proving Theorem~\ref{canG} we saw the following.
\begin{corollary}\label{degmu}
  The curvature of an admissible log function on a line bundle $\L$ is
  $\deg(\L) \cdot \mu$.
\end{corollary}

\section{The Faltings volume on the determinant of cohomology}
\label{sec:volume}

In classical Arakelov theory one defines a volume on the determinant
of cohomology of a line bundle (or, more generally, of a vector
bundle). This data then enters into the Riemann-Roch
theorem. In~\cite{Fal84} Faltings constructs the volume for line bundles in
an axiomatic way. It is also possible to obtain a volume using
analysis (analytic torsion). Here we follow the approach of
Faltings. It would be very interesting if one could
also find a definition using an analogue of analytic torsion but we
have no idea how to do this.

Let $X$ be a complete non-singular curve over $\qpc$. We fix a Green
function $G$ on $X$ out of the almost canonical class. Let $\L$ be a
line bundle on $X$. Recall that the determinant of cohomology of $\L$
is given by
\begin{equation*}
  \lambda(\L)= \det(H^0(X,\L))\otimes \det(H^1(X,\L))^{\otimes -1}\;,
\end{equation*}
where $\det$ is the top exterior power.
\begin{proposition}\label{detcond}
  There exist a correspondence
  \begin{equation*}
    (\L,\llog) \mapsto \text{ log function }
    \log_{\lambda(\L)}^{(\llog)} \text{ on } \lambda(\L)
  \end{equation*}
  from line bundles with an admissible log function with
  respect to $G$ to ``metrized lines'',
  such that the following properties are satisfied:
  \begin{enumerate}
  \item an isometry $\L\to \L'$ induces an isometry $\lambda(\L) \to
    \lambda(\L')$.
  \item \label{detcond2} The behavior with respect to scaling is such that
    \begin{equation*}
      \log_{\lambda(\L)}^{(\llog+\alpha)}=\log_{\lambda(\L)}^{(\llog)}
    + \chi(\L) \cdot \alpha\;,
  \end{equation*}
    where $\chi(\L)$ is the Euler characteristic of $\L$.
  \item \label{detcond3} The canonical isomorphism
    \begin{equation*}
      \lambda(\O(D))\isom \lambda(\O(D-P))\otimes \O(D)[P]\;, 
    \end{equation*}
    where $\O(D)[P]$ is the fiber of $\O(D)$ at $P$, is an isomorphism.
  \end{enumerate}
  Furthermore, these properties determine
  $\log_{\lambda(\L)}^{(\llog)}$ up to common scaling for all $\L$
  together.
\end{proposition}
\begin{proof}
The proof proceeds in a similar manner to the corresponding proof
in~\cite{Fal84}. The uniqueness is clear since we can pass from one
metrized line
bundle to any other by either adding or deleting points or scaling,
so fixing $\log_{\lambda(\L)}^{(\llog)}$ for one metrized line
bundle $\L$ determines it on all of them. Following Faltings again,
fixing a divisor $E$ of degree $r+g-1$, for sufficiently large $r$,
and fixing the log function on the determinant of
$\O(D)$ we get a line bundle $\NN$ on $X^r$ whose fiber at
$(P_1,\ldots,P_r)$ is $\lambda(\O(E-\sum P_i))$. The line bundle
$\NN$ carries a
pseudo-log function determined by our conditions. As shown by
Faltings, $\NN$ is the pullback from $\Pic_{g-1}(X)$, the jacobian
variety of line bundles of degree $g-1$ on $X$, of
$\O(-\Theta)$ under the map $\varphi$ sending $(P_1,\ldots,P_r)$ to $E-\sum
P_i$, where $\Theta$ is the theta divisor of line bundles with a
global section.  It suffices to show that the quasi-log
function on $\NN$ is the pullback of a log function on
$\O(-\Theta)$.

The forms and cohomology classes $\omega_i$ and $\bo_i$ are pullbacks
of classes, which we denote by the same notation, on $\Pic_{g-1}(X)$.
It is known that
\begin{equation*}
  ch_1(\O(-\Theta))= -\sum_{i=1}^g \bo_i \cup \omega_i\;.
\end{equation*}
It follows from Proposition~\ref{excurve} that on $\O(-\Theta)$ there
exist a log function whose curvature form is
\begin{equation*}
    \Curve (\O(-\Theta))= -\sum_{i=1}^g \bo_i \otimes \omega_i\;.
\end{equation*}
Let $p_k: X^r \to X$, $1\le k\le r$ be the projection on the $k$th
factor and let $p_{km}$ be the projection on the $k$ and $m$
factors. Then the curvature on the pullback log function on
$\NN$ is
\begin{equation*}
  - \sum_{i=1}^g \left(\sum_{k=1}^r p_k^\ast (\bo_i)\right) \otimes
  \left(\sum_{k=1}^r p_k^\ast (\omega_i)\right)\;.
\end{equation*}
We now show that the
pseudo-log function on $\NN$ is indeed a log function and
compute its curvature. Let $\NN_m$, for $0\le m\le r$, be the
line bundle whose fiber at $(P_1,\ldots,P_r)$ is $\lambda(\O(E-\sum_{i=1}^m
P_i))$. Then $\NN_0$ is the constant line bundle
$\lambda(\O(E))$ and condition \eqref{detcond3} implies an isometry
$\NN_{m}=\NN_{m-1}\otimes \L_m^{-1}$, where $\L_m$ is the line bundle
whose fiber at $(P_1,\ldots,P_r)$ is the fiber at $P_m$ of
$\O(E-\sum_{i=1}^{m-1} P_i)$. We have
\begin{equation*}
  \L_m = p_m^\ast \O(E) \otimes \bigotimes_{k=1}^{m-1} p_{km}^\ast
  \O(\Delta)^{-1}\;.
\end{equation*}
Thus we obtain an isometry
\begin{equation*}
  \NN\isom \lambda(\O(E))\otimes \bigotimes_{m=1}^r  p_m^\ast
  \O(E)^{-1} \otimes\bigotimes_{k<m} p_{km}^\ast
  \O(\Delta)\;.
\end{equation*}
In particular, the pseudo-log function on $\NN$ is indeed a log
function and we may also compute its curvature, using
Corollary~\ref{degmu} for the
curvature of $\O(E)$, to be (here we follow~\cite[p. 146]{Lan88})
\begin{align*}
  \Curve(\NN)&= - (r+g-1) \sum_{m=1}^r  p_m^\ast \mu + \sum_{k<m}
  p_{km}^\ast \Phi \\
    &= - (r+g-1) \sum_{m=1}^r  p_m^\ast \mu \\
    &\phantom{=} +\sum_{k<m} \left( p_k^\ast \mu + p_m^\ast \mu
     - \sum_{i=1}^g \left(p_k^\ast \bo_i \otimes  p_m^\ast\omega_i
      +p_m^\ast \bo_i \otimes  p_k^\ast\omega_i
    \right)\right)\;.\\
\intertext{The term $(r-1)  \sum_{m=1}^r  p_m^\ast \mu$ cancels,
    leaving us with}
  &= - \sum_{m=1}^r  \sum_{i=1}^g p_m^\ast (\bo_i \otimes \omega_i) -
  \sum_{k<m} 
  \sum_{i=1}^g \left(p_k^\ast \bo_i \otimes  p_m^\ast\omega_i
      +p_m^\ast \bo_i \otimes  p_k^\ast\omega_i
    \right) \\
    &=  - \sum_{i=1}^g \left(\sum_{k=1}^r p_k^\ast (\bo_i)\right) \otimes
  \left(\sum_{k=1}^r p_k^\ast (\omega_i)\right)
\end{align*}

Thus, $\log_{\NN}$ and $\varphi^\ast \log_{\O(-\Theta)}$ have the
same curvature hence they differ by the integral of a holomorphic
form $\omega$. However, both log functions are invariant with respect to the
action of the symmetric group on $X^r$. It follows that $\omega =
\sum_i p_i^\ast \omega'$ where $\omega'\in \Omega^1(X)$. But then
$\omega$ can be pulled back from $\Pic_{g-1}$ and therefore
$\log_{\NN}$ can be pulled back as well, which proved the
result.
\end{proof}

\section{Relations with the theory of Coleman and Gross}
\label{sec:colgros}

In~\cite{Col-Gro89} Coleman and Gross define a $p$-adic height pairing on
curves with good reduction above $p$ as a sum of local terms.
In \cite{Bes99a} we prove that this local height pairing is the same as the
one defined by \nekovar\ in~\cite{Nek93}. The
height pairing is defined for divisors of degree $0$ and one expects
that it coincides with the restriction to these divisors of the
Arakelov intersection pairing. Once we define the Arakelov
intersection we will prove that this is indeed the case. At the moment
we can only prove that the local terms above $p$ agree.

We first recall the local theory in~\cite{Col-Gro89}. We reformulate
slightly since Coleman 
and Gross work over $\C_p$, but we can easily work over $\qpc$
instead. Also, their definition is only for curves with good reduction
but once one has integration theory in the bad reduction case as well
the extension is done verbatim. Let $X$ be again a complete non
singular curve over $\qpc$. Recall from Section~\ref{sec:double} the
space $\T$ of forms 
of the third kind on $X$, the subspace $\T_l$ of dlog forms and the
map $\lag: \T/\T_l \to \hdr^1(X)$.

The theory of Coleman and Gross depends, as does our theory, on the
choice of a subspace 
$W\in \hdr^1(X)$ complementary to $\Omega^1(X)$ which is isotropic with
respect to the cup product. In~\cite{Col-Gro89} this is not needed but
is imposed if one wants to make the height
pairing symmetric.
\begin{definition}
  For any divisor $D$ of degree $0$ on $X$  we let
  $\omega_D\in 
  \T$ be the unique form satisfying $\res(\omega_D)=D$ and
  $\lag(\omega_D)\in W$.
\end{definition}
The uniqueness of $\omega_D$ follows since a form in $\T$ with zero
residue divisor is holomorphic, and on such a form $\lag$ is the identity.
\begin{definition}
  The local height pairing of Coleman and Gross is defined as follows:
  Let $D$ and $E$ be two divisors of degree $0$ on $X$ with disjoint
  supports. Then their
  pairing is given by $\pair{D,E}:=\int_{E}\omega_D$. The integral in
  the definition is the Coleman integral of $\omega_D$ evaluated in
  the standard way on $E$.
\end{definition}

It is now clear that the equality of the local height pairing and the
Arakelov pairing at primes above $p$ follows from the following result.
\begin{theorem}\label{sameasColGro}
  Let the space $W$ be chosen. Then for any divisor $D$ of degree $0$
  on $X$ we have $\dd  G_D =\omega_D$.
\end{theorem}
\begin{proof}
We have $\dbd (G_D)|_{X-D}=0$ by Corollary~\ref{degmu} so that $\dd  G_D$
is holomorphic
outside $D$. It follows from the definition of $G_D$ that it has
logarithmic singularities and its residue divisor is exactly $D$.
Let $\omega_D^\prime:= \dd  G_D -\omega_D$. Then $\omega_D^\prime\in
\Omega^1(X)$ for each $D$. Further, by Proposition~\ref{chidef1} we
have $\omega_D^\prime=0$
for a principal divisor $D$, so the map $D\mapsto \omega_D^\prime$
factors through $J$, where $J$ is the Jacobian of $X$, and it is
clearly additive. To prove that
$\omega_D^\prime=0$ it suffices to show that for any $w\in W$ we have
$w\cup \omega_D^\prime=0$. 
The projection $\lag$ is the identity on $\Omega^1(X)$ and
by construction it maps $\omega_D$ to $W$. Since $W$ is isotropic we
have $w\cup \omega_D^\prime= w\cup \lag(\dd  G_D)$, and by
Corollary~\ref{doubleisproj} this equals $\pair{w,\dd  G_D}_\gl$.
The map $D\mapsto \pair{w,\dd  G_D}_\gl$ is
an additive map on $J$, so it
suffices to prove that it is locally analytic and its derivative is
$0$. for any $w\in W$. We consider the map 
$X^r \to J$, $(P_1,\ldots P_r)\mapsto \sum P_i -r P_0$, for some
$P_0$, which is surjective for sufficiently large $r$. It will suffice
to show that the map $(P_1,\ldots P_r) \mapsto \pair{w,\dd  G_{\sum P_i-r
  P_0}}_\gl$ has zero derivative with respect to every $P_i$, and for
this it suffices to check the derivative of $P\mapsto \pair{w,\dd 
  G_{P-P_0}}_\gl$. Let
$\partial/\partial P$ be a vector field on $X$. By Lemma~\ref{derdoub}
we have
\begin{equation*}
  \frac{\partial}{\partial P} \pair{w,\dd  G_{P-P_0}}_\gl =
  \pair{w, \frac{\partial}{\partial P} \dd  G_{P-P_0}}_\gl=
  \pair{w, \frac{\partial}{\partial P} \dd  G_{P}}_\gl\;.
\end{equation*}
We will now determine the cohomology class of the form
$(\partial/\partial P) \dd  G_{P}$. The situation can be described as
follows. We have on $X\times X-\Delta$ the Coleman form $\dd  G(P,Q)$. We
differentiate with respect to the vector field $\partial/\partial P$
on the first variable $P$
and then restrict to the fiber at $P$. Since $\dd  G$ is closed we have
by \eqref{yek} $(\partial/\partial P) \dd  G = \dd (\dd
G|_{\partial/\partial P})$. Since the
retraction operator $|_{\partial/\partial P}$ is $\O_X$-linear we have
by Lemma~\ref{dercurve} that $\db (\dd  G|_{\partial/\partial P})= (\dbd
G)|_{\partial/\partial P}$, where the retraction on $\hdr^1(X\times
X-\Delta) \otimes \Omega^1(X\times X-\Delta)$ operates on the second
factor. By Corollary~\ref{dbdG} and Definition~\ref{cancurve} we have
\begin{equation*}
  \db (\dd  G|_{\partial/\partial P})= \Phi|_{\partial/\partial P}=
  \frac{1}{g} \sum_{i=1}^g \pi_1^\ast \bo_i \otimes
  \pi_1^\ast(\omega_i|_{\partial/\partial P})
    - \sum_{i=1}^g \pi_2^\ast \bo_i \otimes
    \pi_1^\ast(\omega_i|_{\partial/\partial P})
\;,
\end{equation*}
which, when restricted to the fiber above $P$ yields $\sum 
\bo_i \otimes \alpha_i$, with $\alpha_i$ the constant
$-(\omega_i|_{\partial/\partial P})(P)$. If we now
represent the $\bo_i$ by forms of the second kind and restrict further
to $U\subset X$ where all these forms are holomorphic, then it follows
from Proposition~\ref{volddbsurjects} and the description of the
$\db$ operator on affine varieties that
\begin{equation*}
   ((\dd  G)|_{\partial/\partial P})|_{P\times U}= \sum \alpha_i \int
  \bo_i +f
\end{equation*}
where $f\in \O(U)$. Therefore,
\begin{equation*}
   \dd ((\dd  G)|_{\partial/\partial P})|_{P\times U}= \sum \alpha_i 
  \bo_i +\dd  f\;,
\end{equation*}
whose cohomology class belongs to $W$ by the definition of the $\bo_i$.
It follows that $(\partial/\partial P) \dd  G_{P}$ is a form of the
second kind 
representing a class in $W$ and the isotropy of $W$ completes the proof.
\end{proof}
The following proposition demonstrates that the Green function is
forced on us if we assume compatibility with the Coleman-Gross
pairing, symmetry and a natural residue condition, analogous to the one
we have in the classical theory. As mentioned in the introduction,
this was our original approach. At the same time we provide a formula
for the Green function using only Coleman integration in one variable.
\begin{proposition}\label{charofgr}
  The canonical Green function $G$ is the unique function up to
  constant satisfying the following properties.
  \begin{enumerate}
  \item $G$ is symmetric
  \item The induced height pairing \eqref{heightpairing} is the
    Coleman-Gross height pairing.
  \item \label{rescond} The following residue condition is satisfied:
    For any point 
    $P$ the canonical map $\bfo_X \otimes \O(P) \to (\qpc)_P$, where $
    (\qpc)_P$ is the skyscraper sheaf at $P$ with fiber $\qpc$, 
    given by $\omega\otimes f\to \res_P(f\omega)$, is an isometry.
  \end{enumerate}
\end{proposition}
\begin{proof}
That $G$ satisfies the 3 conditions follows from
Theorem~\ref{sameasColGro} and Theorem~\ref{canG}. For uniqueness
we will in fact prove the following explicit formula for the Green
function in terms of the Coleman-Gross height pairing:
Choose $a$ and $b$ two
points in $X$. Then we have
\begin{equation}
  \label{Green-formula}
  G(P,Q)=\frac{1}{2g}\left(
    \int_{2gP-\Div \omega_2-Q-b}\omega_{Q-b}+
    \int_{2gb-\Div \omega_1-P-a}\omega_{P-a}
  \right)
\end{equation}
where $\omega_1$ (respectively $\omega_2$) is any form with log
singularities at $P$ and $a$ (respectively $Q$ and $b$) and such
that the log of its residues at these two points is the same.

Let $P$ be a point of $X$ and let $f$ be a local parameter at $P$.
Let $\omega$ be a form with a simple pole at $P$. We can write
$\omega=(1/f) f\omega$ and by condition \eqref{rescond} we should have
\begin{align*}
\log (\res_P(\omega))&=(\log_{\O(P)}(1/f) + \log_{\bfo_X}(f\omega))_P
\\&=(G_P+\log(f) -\log(f)+ \log_{\bfo_X}(\omega))_P\\&=\lim_{z\to P}
G_P(z)+\log_{\bfo_X}(\omega)(z)\;.
\end{align*}
Since the log function on $\bfo$ is
admissible we have
\begin{equation*}
  \log (\res_P(\omega))=\lim_{z\to P} G_P(z) + G_{\Div(\omega)}(z)
  +\ilog_\bfo(\omega)=G_{\Div(\omega)+P}(P)+\ilog_\bfo(\omega)\;.
\end{equation*}
Consider now a any two points $P$ and $Q$ on $X$ and a
differential $\omega$ with simple poles at both $P$ and $Q$ such that
the logs
of the residues of $\omega$ at $P$ and $Q$ are equal. We find the
equations
\begin{equation*}
  G_{\Div(\omega)+P+Q}(P)-G_Q(P)+\ilog_\bfo(\omega)=
  G_{\Div(\omega)+P+Q}(Q)-G_P(Q)+\ilog_\bfo(\omega)
\end{equation*}
and from the symmetry of $G$ we get by subtracting
\begin{equation*}
  G_{\Div(\omega)+P+Q}(P-Q)=0\;,
\end{equation*}
from which we get, again by symmetry
\begin{equation*}
  G_{P-Q}(\Div(\omega)+P+Q)=0.
\end{equation*}
Note that the divisor $\Div(\omega)+P+Q$ has degree $2g$. It follows
that for any point $x$ in $X$ different from $P$ and $Q$ we have
\begin{align*}
  2g G_{P-Q}(x)&=G_{P-Q}(2gx)=G_{P-Q}(2gx-(\Div(\omega)+P+Q))\\&=
  \int_{2gx-(\Div(\omega)+P+Q)}\omega_{P-Q}\;.
\end{align*}
This gives us formula \eqref{Green-formula} as follows: Write the
bilinear height pairing as $\pair{\bullet,\bullet}$. We choose a
constant for our Green function by insisting that $G(a,b)=0$. Then
\begin{align*}
  2g G(P,Q) &= 2g\pair{P,Q} = 2g (\pair{P,Q-b}+\pair{b,P-a}+\pair{b,a})\\&=
      \int_{2gP-\Div \omega_2-Q-b}\omega_{Q-b}+
      \int_{2gb-\Div \omega_1-P-a}\omega_{P-a}+0\;.
\end{align*}
\end{proof}

\section{Local theory over non algebraically closed fields}
\label{sec:lonnalg}

So far we found it more convenient to develop the local theory over
$\qpc$. However, for Arakelov theory we will need to work with finite
extensions of $\Q_p$. In this section we collect all the necessary
results needed for doing this.

Suppose now that $K$ is a finite extension of $\Q_p$ and that $X/K$ is
a smooth complete curve of genus $\ge 1$. We assume that the space $W$
is also defined over $K$. Finally, we choose a branch of the logarithm
defined over $K$.
\begin{proposition}
  There exist a canonical Green function $G$ for $X$ defined over $K$. It
  is defined up to a constant in $K$.
\end{proposition}
\begin{proof}
The form $\dd G$ is uniquely determined by the conditions spelled
out in Theorem~\ref{canG}. If $\sigma\in \Gal(\overline{K}/K)$, then
all of these conditions are invariant under $\sigma$, so
$\sigma(\dd G)=\dd G$. By Proposition~\ref{2.9} the form $\dd G$ is in fact
defined over
$K$, hence it has an integral defined over $K$ and defined up to a
constant in $K$.
\end{proof}

From now on we will assume that a Green function defined over $K$ has
been fixed.
\begin{corollary}
  If $\L$ is a line bundle on $X$, then there exists an
  admissible log function on $\L$ defined over $K$. If $D$ is a
  divisor on $X$, then the canonical log function on $\O(D)$ of
  Definition~\ref{odcan} is defined over $K$.
\end{corollary}

Using our Green function we can define the local intersection pairing.
\begin{definition}\label{locpair}
  Let $D$ and $E$ be divisors on $X$ with disjoint supports and let
  $\overline{D}=\sum n_i P_i$ and
  $\overline{E}=\sum m_j Q_j$ be their extensions to
  $\overline{X}$. The \emph{local intersection pairing} $\pair{D,E}\in K$ is
  defined by 
  \begin{equation*}
    \pair{D,E}=\sum n_i m_j G(P_i,Q_j)\;.
  \end{equation*}
\end{definition}
Hidden in this definition is the fact that the pairing indeed takes
values in $K$, which follows trivially from the properties of $G$.

Until the end of this section, we develop the relation between log
functions and determinants. The following definition of log functions
is a mere specialization of previous definitions for the case of
dimension $0$.

\begin{definition}
  Let $K'$ be a finite extension of
  $K$. Let $V$ be a $K'$-line, i.e., a one-dimensional $K'$ vector
  space. A log function on $V$ gives, for any embedding
  $\tau: K' \inject \overline{K}$ fixing $K$ a log function $\log_\tau$
  on $V_\tau:=V\otimes_{K'} \overline{K}$. Such a log function is said
  to be defined over $K$ if for any
  automorphism $\sigma\in\Gal(\overline{K}/K)$ we have $\sigma (\log_\tau(x)) =
  \log_{\sigma\tau}(\sigma(x))$ for $x\in V_\tau$, where $\sigma: V_\tau \to
  V_{\sigma\tau}$ is the evident map.
\end{definition}
Note that a log function on $V$ defined over $K'$ is simply a log
function $\log_V:V\to K'$. 
Note also that if $\L$ is a metrized line bundle over a $K$-variety $X$ and
$\llog$ is a log function on $\L$. Then the fiber of $\L$ over any
closed points acquires in a natural way a log function over $K$.

We will now consider a $K'$-line $V$ as above. Since $V$ is a finite
dimensional vector space over $K$, its determinant $\det V=\det_K V$
is defined. We want to obtain log functions on the determinants in
certain situations.
\begin{definition}\label{detlog}
  Suppose that $U$ and $V$ are both $K'$-lines with log functions
  defined over $K'$. Then, the
  $K$-line $(U:V):=\det(U)\otimes \det(V)^{-1}$ has a log function
  defined over $K$ in
  the following way: Let $\alpha: U\to V$ be an isomorphism such that
  $\log_V\circ \alpha=\log_U+c$, where $c\in K'$. The isomorphism
  $\alpha$ induces by composition a canonical isomorphism $\beta: (U:V) \to K$
  and we define $\log_{(U:V)}= \log\circ \beta - \tr_{K'/K} c$.
\end{definition}
The log function we defined is independent of the choices made. Indeed
if $\alpha$ is changed to $c' \alpha$, then $c$ is changed to
$c+\log(c')$ while $\beta$ is changed to $(\norm_{K'/K} c') \beta$ so
$\log_{(U:V)}$ is unchanged. Another way of describing this log
function is to say that if $\alpha$ is an isometry $\beta$ is also is
an isometry while if we scale the log function on $U$ by $c$ then we
scale the log function on $(U:V)$ by $\tr_{K'/K} c$.

Now we consider extension of scalars.
\begin{definition}\label{logscext}
  Suppose $V$ is a $K'$-metrized line and $K''$ is a finite extension
  of $K'$. We define the log function on the 
  $K''$-line $W:=V\otimes_{K'} K''$
  by $\log(v\otimes \alpha)= 
  \log(v)+\log(\alpha)$.
\end{definition}

Determinants behave in a well known way under extensions. Suppose
$[K'':K']=n$. We have a canonical isomorphism
\begin{equation}
  \label{blubird}
  \det_K(W)\otimes \det_K(K'')^{-1}\isom (\det_K(V)\otimes
  \det_K(K'))^{\otimes n}
\end{equation}
defined as follows: Choose a $K'$-isomorphism $K'^{\oplus n}
\xrightarrow{\sim} K''$. This induces an isomorphism $V^{\oplus n}
\xrightarrow{\sim} V\otimes_{K'} K''$ and as consequence isomorphisms
\begin{align*}
  (\det_K K')^{\otimes n} &\xrightarrow{\sim} \det_K K''\\
  (\det_K V)^{\otimes n} &\xrightarrow{\sim} \det_K (V\otimes_{K'} K'')\;,
\end{align*}
from which \eqref{blubird} follows. It is easily verified that this
isomorphisms is independent of the choices. The behavior with respect
to log functions is also easily checked.
\begin{proposition}\label{logscext1}
  In the situation described above the canonical
  isomorphism \eqref{blubird} is an isometry.
\end{proposition}
\begin{proof}
This is clear from the description above if we choose an isometry $K'
\xrightarrow{\sim} V$. Observing how the log functions change with
respect to scaling finishes the proof.
\end{proof}

For log functions defined over $K$ we can define the log function on
the determinant, and not only on the quotient of two determinants.
\begin{lemma}\label{goodlog}
  Suppose $K'$ is a finite extension of $K$, $V$ a $K'$-line equipped
  with a log function over $K$. Then the $K$-line $\det_{K}V$ has a
  unique log function satisfying the following property: The
  isomorphism
  \begin{equation*}
    V\otimes_K \bar{K} \isom \sum_{\sigma:K' \to \bar{K}}
    V\otimes_{K',\sigma} \bar{K}\;,\quad x\otimes \alpha \mapsto \sum
    x\otimes \alpha\;,
  \end{equation*}
  where the sum is over all embeddings $\sigma:K' \to \bar{K}$ fixing
  $K$, induces an isomorphism
  \begin{equation*}
    (\det_{K}V)\otimes_K \bar{K} \isom \det_{\bar{K}} V\otimes_K
    \bar{K} \isom  \otimes_{\sigma:K' \to \bar{K}}
    V\otimes_{K',\sigma} \bar{K}
  \end{equation*}
  and this isomorphism becomes an isometry with respect to the log
  functions on both sides.
\end{lemma}
\begin{proof}
Uniqueness is clear. To prove existence, let $x$ be a basis of $V$
over $K'$ and let $\{\beta_i\}$, $i=1,\ldots n$ be a basis of $K'$ over
$K$. Then $\{ \beta_i x\}$ is a basis of $V$ over $K$. Number the
embeddings $\sigma_j$, $j=1,\ldots, n$. We have
\begin{equation*}
  \beta_i x\mapsto \beta_i x\otimes 1 \mapsto \oplus_j \beta_i
  x\otimes 1 = \oplus_j \sigma_j(\beta_i)\cdot  (x\otimes 1)_j 
\end{equation*}
where the subscript $j$ is given to distinguish the different
components. The basis $\wedge_i (\beta_i x)$ of $\det_K V$ is mapped
to $\det(\sigma_j(\beta_i))\cdot \wedge_j (x\otimes 1)_j$ and this
forces us to define
\begin{equation*}
  \log(\wedge_i (\beta_i x))= \log (\det(\sigma_j(\beta_i)))+\sum_j
  \log_j((x\otimes 1)_j)\;,
\end{equation*}
where $\log_j$ is the log function on
$V\otimes_{K',\sigma_j}\bar{K}$. Since the log function is defined
over $K$ we have $\log_j((x\otimes 1)_j)=\sigma_j(\log(x))$ so we
obtain
\begin{equation*}
  \log(\wedge_i (\beta_i x))= \log
  (\det(\sigma_j(\beta_i)))+\tr_{K'/K}\log(x)\;.
\end{equation*}
The proof will be complete if we show that the $\log
(\det(\sigma_j(\beta_i)))\in K$. This is clear since applying an
automorphism of $\bar{K}$ over $K$ multiplies
$\det(\sigma_j(\beta_i))$ by $\pm 1$ so the log is unchanged.
\end{proof}
The log function just defined is easily seen to be compatible with the
one defined in Definition~\ref{detlog} as follows:
\begin{proposition}\label{9.13}
  If $U$ and $V$ are two $K'$-lines with a log function defined over
  $K$. Let $(U:V)$ be the metrized $K$-line of
  Definition~\ref{detlog}. Then the isomorphism $(U:V)=\det(U)\otimes
  \det(V)^{-1}$ is an isometry, with $\det(U)$ and $\det(V)$ having
  the log functions defined in Lemma~\ref{goodlog}.
\end{proposition}
\begin{proposition}\label{9.14}
  Let $K'$ be a finite extension of $K$. Let $\log$ be the log
  function on $K'$ obtained from the one on $K$. Let $V=K'$ with this
  log function (defined over $K$). The trace form induces an
  isomorphism $\det_K(V) \otimes \det_K(V) \to K$ and this is an isometry.
\end{proposition}
\begin{proof}
Let $\{\beta_i\}$ be a basis of $K'$ over $K$ and let $\beta_i^\prime$
be a dual basis with respect to the trace form. The log function we
defined sends $\wedge_i \beta_i$ to $\log(\det(\sigma_j(\beta_i)))$,
where $\sigma_j$ are the embeddings of $K'$ in $\bar{K}$,
and similarly with $\beta_i^\prime$ replacing $\beta_i$. The duality
with respect to the trace form implies $(\sigma_j(\beta_i)) \cdot
(\sigma_j(\beta_i^\prime))^t = I$, hence $\log (\wedge_i \beta_i)+\log
(\wedge_i \beta_i^\prime)=0$, which is what we want.
\end{proof}

\section{The intersection pairing}
\label{sec:global}

We now combine the $p$-adic analysis of the previous chapters to
obtain a $p$-adic Arakelov intersection pairing. For motivation to the
setup introduced here the reader is encourages to look
at~\cite{Col-Gro89}.

The general setup is as follows: $F$ is a number field and $p$ is a
prime. We choose a ``global log''- a continuous idele class character
\begin{equation*}
  \ell: \mathbb{A}_F^\times/F^\times \to \Q_p\;.
\end{equation*}
One deduces from $\ell$ the following data:
\begin{itemize}
\item For any place $v\ndiv p$ we have $\ell_v(\O_{F_v}^\times)=0$ for
  continuity reasons, which implies that $\ell_v$ is completely
  determined by the number $\ell_v(\pi_v)$, where $\pi_v$ is any
  uniformizer in $F_v$.
\item For any place $v|p$ one can decompose
  \begin{equation*}
    \xymatrix{
    {\O_{F_v}^\times}  \ar[rr]^{\ell_v} \ar[dr]^{\log_v} & &   \Q_p\\
        & F_v\ar[ur]^{t_v}
          }
  \end{equation*}
  where $t_v$ is a $\Q_p$-linear map. We assume that $\ell_v$ is
  ramified in the 
  sense that it does not vanish on $\O_{F_v}^\times$.
  It is then possible to extend $\log_v$ to $F_v^\times$ in such a way
  that the diagram above remains commutative when $\O_{F_v}^\times$ is
  replaced by $F_v^\times$, i.e., we have the decomposition
  \begin{equation}
    \label{ldecomp}
    \ell_v = t_v \circ \log_v\;.
  \end{equation}
\end{itemize}
Note that for $v\ndiv p$ and for any $\O_{F_v}$-ideal $I$ in $F_v$ we
can define unambiguously
\begin{equation*}
  \ell_v(I):=\ell_v(\tau),\quad \tau \text{ is a generator of } I\;.
\end{equation*}
Let $\XX$ be an arithmetic surface over $\O_F$ (i.e., a proper regular curve
over $\O_F$). Let $X_v:= \XX \otimes_{\O_F} F_v$. We make the following
additional choices for each $v|p$
\begin{itemize}
\item A space $W_v$ in $\hdr^1(X_v/F_v)$ complementary to $F^1$ as in
  Section~\ref{sec:canonical}.
\item A choice of a Green function $G_v$, defined over $F_v$, out of
  the almost canonical one (choice of a constant).
\end{itemize}

\begin{definition}
  A \emph{$p$-adic Arakelov divisor} (Arakelov divisor for short) on $\XX$ is
  a formal combination 
  \begin{equation*}
    D=D_{\fin}+D_\infty,\quad \text{where } D_\infty=\sum_{v|p}
    \lambda_v X_v\;,
  \end{equation*}
  where  $\lambda_v\in F_v$. Here, $X_v$ should be treated as a formal
  symbol. The group of all Arakelov divisors on $\XX$ is
  denoted $\divar(\XX)$.
\end{definition}
\begin{definition}
  Let $D$ and $E$ be two Arakelov divisors and suppose that the
  intersections of $D_\fin$ and $E_\fin$ with the generic fiber have
  disjoint supports. The \emph{Arakelov intersection pairing} of $D$ and $E$
  is defined as
  \begin{equation*}
    D \cdot E =\sum_v [D,E]_v
  \end{equation*}
  where the local intersection multiplicities $[D,E]_v\in \Q_p$ are
  defined by the following rules (extended by symmetry):
  \begin{enumerate}
  \item If $v\ndiv p$, then
    \begin{equation*}
      [D,E]_v = \ell_v(\pi_v) \pair{D_\fin,E_{\fin}}_v
    \end{equation*}
    where $\pair{D_{\fin},E_{\fin}}_v$ is the usual intersection
    multiplicity at $v$ of the finite parts of $D$ and $E$.
  \item If $v|p$, then we have
    \begin{equation*}
      [D,E]_v=t_v(\pair{D,E}_v)
    \end{equation*}
    where the intersection multiplicities $\pair{D,E}_v\in F_v$ are given
    by the following rules:
    \begin{enumerate}
    \item if $w\ne v$, then $\pair{D,\lambda_w X_w}_v=0$.
    \item $\pair{\lambda_1 X_v, \lambda_2 X_v}_v=0$.
    \item if $D$ is a finite divisor, then $\pair{D,\lambda X_v}_v =
      \lambda \deg D_F$, where $D_F$ is the generic part of $D$.
    \item Suppose $D$ and $E$ are finite and let $D_v$
      and $E_v$ be
      their images in $X_v$. Then we have
      \begin{equation*}
        \pair{D,E}_v=\pair{D_v,E_v}\;,
      \end{equation*}
      where this last pairing is the one of Definition~\ref{locpair}
      taken with respect to the Green function at $v$.
    \end{enumerate}
  \end{enumerate}
\end{definition}
\begin{definition}
  An Arakelov line bundle on $\XX$ is a line bundle $\L$ on $\XX$
  together with a choice of an admissible metric on $\L_v$ for every
  $v|p$. The trivial Arakelov line bundle $\O_{\XX}$ is the line bundle
  $\O_{\XX}$ together with the canonical metric.
\end{definition}
There is are obvious notions of isomorphisms of Arakelov line bundles
and of the tensor product of them.
\begin{definition}
  Let $\L$ be an Arakelov line bundle on $\XX$ and let $s$ be a
  rational section of $\L$. The Arakelov divisor $(s)$ of $s$ is
  defined as $(s)=(s)_\fin+ (s)_\infty$ where $(s)_\fin$ is the usual
  divisor of $s$ and
  \begin{equation*}
    (s)_\infty = \sum_{v|p} \ilog_\L(s_v) X_v\;.
  \end{equation*}
  In particular, considering the case $\L=\O_{\XX}$ we obtain the
  Arakelov divisor of a rational function.
\end{definition}
Clearly we have
\begin{equation*}
  (s\otimes t)=(s)+(t)
\end{equation*}
for sections of two line bundles. In particular we have $(fg)=(f)+(g)$
for any two functions and $(fs)=(f)+(s)$ where $f$ is a rational
function and $s$ a section of a line bundle.
\begin{definition}
  The group of principal Arakelov divisors is the group
  \begin{equation*}
    \prinar(\XX) := \{(f)|\; f \in F(\XX)^\times\}\;.
  \end{equation*}
\end{definition}
\begin{lemma}
  We Let $D$ be a finite divisor and suppose $v\ndiv p$. Then
  $[D,(f)]_v=\ell_v(f(D))$.
\end{lemma}
\begin{proof}
Well known, see for example in~\cite{Col-Gro89} Proposition~1.2 and
its proof.
\end{proof}
\begin{proposition}
  If $f$ is a rational function and $D$ and Arakelov divisor, then
  \begin{equation*}
    D\cdot (f)=0
  \end{equation*}
\end{proposition}
\begin{proof}
The only interesting case is when $D$ is finite, where we have
\begin{align*}
  D\cdot (f)&=\sum_{v\ndiv p}[D,(f)]_v + \sum_{v|p}\left(
    [D,(f)_\fin]_v + [D,\ilog (f_v) X_v]_v\right)\\
   &= \sum_{v\ndiv p} \ell_v(f(D))+ \sum_{v|p} t_v\left(
    G_{(f_v)}(D) + \ilog (f_v) \deg D_F\right)\\
   &= \sum_{v\ndiv p} \ell_v(f(D))+ \sum_{v|p} t_v (\log_v(f(D)))
   \quad \text{by Definition~\ref{ilogdef}}\\
   &= \sum_{v\ndiv p} \ell_v(f(D))+ \sum_{v|p} \ell_v(f(D))
   \quad \text{by \eqref{ldecomp}}\\
   &= 0
\end{align*}
since $\ell$ is an idele class character.
\end{proof}
\begin{definition}
  The \emph{Arakelov Chow group} is the quotient group 
  \begin{equation*}
    \chowar(\XX):=\divar(\XX)/\prinar(\XX)\;.
  \end{equation*}
\end{definition}
The following result is now standard
\begin{proposition}
  There is a unique bilinear Arakelov intersection pairing on
  $\chowar(\XX)$ specializing to the previously defined intersection
  pairing for two divisors with disjoint supports on the generic fiber.
\end{proposition}
\begin{definition}
  The \emph{Arakelov Picard group} of $\XX$ is the group $\picar(\XX)$ of
  isometry classes of line bundles on $\XX$ with admissible metrics at
  primes above $p$.
\end{definition}
\begin{definition}
  Given an Arakelov divisor $D=D_\fin + \sum \lambda_v X_v$, we define
  the metrized line bundle
  $\O(D)$ on $\XX$ as follows: As a line bundle it is simply
  $\O(D_\fin)$ and if $v|p$, then the log function on
  $\O(D_\fin)_v=\O((D_\fin)_v)$ is the canonical one
  (Definition~\ref{odcan})  scaled by $\lambda_v$
\end{definition}
The line bundle $\O(D)$ is admissible and it is clear that any
admissible metrized line bundle is isomorphic to $\O(D)$ for some
Arakelov divisor $D$.
The following result is clear.
\begin{proposition}
  There is an isomorphism $\picar(\XX) \isom \chowar(\XX) $ given by
  the two inverse maps
  \begin{equation*}
    \L \mapsto c(\L),\quad
    D\mapsto \O(D)\;.
  \end{equation*}
  where
  \begin{equation}
    \label{archern}
    c(\L):= (s),\; s \text{ a rational section of } \L\;.
  \end{equation}
\end{proposition}

\begin{definition}
  Let $\NN$ be a metrized line bundle over $\O_F$, i.e., a
  locally free $\O_F$-module of rank $1$ together with a choice, for
  each $v|p$, of a log function $\log_v$ on $\NN_v:=\NN\otimes_{\O_F} F_v$. We
  define the \emph{degree} of $\NN$ as follows: Fix an isomorphism
  $\theta:F \xrightarrow{\sim} \NN\otimes F$, which induces local
  isomorphisms $\theta_v: F_v \xrightarrow{\sim}
  \NN_v$ for each $v$. Then we define
  \begin{equation*}
    \deg \NN  = \sum_{v|p} t_v(\log_v (\theta_v(1)))-\sum_{v\ndiv p}
    \ell_v (\theta_v^{-1} \NN_v)\;.
  \end{equation*}
\end{definition}
It is very easy to see that this definition is independent of the
choice of the isomorphism $\theta$.

We next generalize the notion of degree to line bundles over finite
$\O_F$-schemes. We use here the theory of the determinant line
bundle~\cite{Knu-Mum76}.
Suppose that $A$ is a finite integral $\O_F$-algebra and that
$\NN$ is a line bundle on $\Spec(A)$. Let $L$ be the fraction field of
$A$. Let $w$ be a place of $L$ above the place $v$ of $F$, lying above
$p$. As before, the choice of $\log_v$ extends uniquely to a branch
$\log_w$ on $L_w$.
\begin{definition}
  We say that $\NN$ is metrized if for any such $w$ we are given a
  log function on $\NN \otimes_A L_w$. We say it is metrized over $F$
  if for each such $w$ lying over $v$ this log function is defined
  over $F_v$.
\end{definition}
\begin{definition}
  Let $\NN$ be a metrized line bundle on $A$. Then, the \emph{degree} of
  $\NN$, $\deg(\NN)$, is defined as the degree of the line bundle
  $\det_{\O_F} \NN \otimes (\det_{\O_F} A)^{-1}$, where the log
  function on
  \begin{equation*}
    (\det_{\O_F} \NN \otimes (\det_{\O_F} A)^{-1})\otimes_{\O_F} F_v
    = \bigotimes_{w|v} \det (\NN\otimes_A L_w) \otimes (\det L_w)^{-1}
  \end{equation*}
  is the tensor product of the log functions of Definition~\ref{detlog}.
\end{definition}
\begin{proposition}\label{9.16}
  We have $\deg (\NN_1 \otimes \NN_2) = \deg (\NN_1) +\deg( \NN_2)$.
\end{proposition}
\begin{proof}
This is clear for line bundles on $\O_F$. For more general line
bundles one can argue as follows: We can find sections $A\to \NN_1$
and $A\to \NN_2$ such that the supports of the cohomology of the
resulting complexes are disjoint (choose the first section arbitrarily and
choose the second to avoid the support of the first). It follows that
the tensor product of the two complexes over $A$ is exact. Taking
determinants we find
\begin{equation*}
  \det \NN_1 \otimes \det \NN_2 \isom \det (\NN_1 \otimes \NN_2)
  \otimes \det A
\end{equation*}
or
\begin{equation*}
  (\det \NN_1\otimes \det(A)^{-1}) \otimes (\det \NN_2 \otimes
  \det(A)^{-1}) \isom \det (\NN_1 \otimes \NN_2) \otimes \det(A)^{-1}\;.
 \end{equation*}
It is easy to see that this isomorphism is an isometry, giving the result.
\end{proof}
\begin{proposition}
  Suppose $f:\Spec(B) \to \Spec(A)$ is a surjective morphism of finite
  $\O_F$-schemes of degree $m$ and $\NN$ is a metrized line bundle on
  $\Spec(A)$. Then $\deg (f^\ast \NN)= m \deg(\NN)$.
\end{proposition}
\begin{proof}
We need to compute
\begin{equation*}
  \det (\NN \otimes_A B) \otimes  \det(B)^{-1} = \det ( B
  \otimes (A\to \NN))^{-1}
\end{equation*}
for every section $A\to \NN$.
Take an injection of $A$ modules $A^m \to  B$ whose cokernel is
supported on a finite number of points. By choosing the section $A\to
\NN$ appropriately, as we did in the proof of Proposition~\ref{9.16},
we can replace $B$ by $A^m$ in the last
equality to get
\begin{equation*}
   \det (\NN \otimes_A B) \otimes  \det(B)^{-1}=
  \det (A^{m} \otimes (A\to \NN))^{-1} = (\det (\NN)\otimes
  \det(A)^{-1})^m\;.
\end{equation*}
It follows immediately from Proposition~\ref{logscext1} that this is
an isometry and the result follows.
\end{proof}
\begin{proposition}
  Let $L/F$ be a finite extension of fields and let $\O_L$ be the ring
  of integers in 
  $L$. Define an idele class character on $L$ by $\ell_L := \ell\circ
  N_{L/K}$. Let $\NN$ be a line bundle on $\Spec(\O_L)$. Then
  $\deg(\NN)$ is the same as the degree of $\NN$ as an $\O_L$ bundle,
  computed with respect to $\ell_L$.
\end{proposition}
\begin{proof}
By applying the previous two results one immediately reduces to the
case where $L/F$ is a Galois extension of degree $m$, $\NN= \M
\otimes_{\O_F} \O_L$ for a line bundle $\L$ on $\O_F$ and one has to
show that the degree of $\NN$, computed with respect to $\ell_L$, is
the same as $m \deg(\M)$. We may choose an isomorphism $\theta: F \to
\M \otimes F$, extend scalars to $L$ and use it to compute the degree
of $\NN$. Let $v$ be a place of $F$ and $w|v$ a place 
of $L$. If $v\ndiv p$ and $\theta_v^{-1} \NN_v = \pi_v^k \O_{F_v}$, then also
$\theta_w^{-1} \NN_w = \pi_v^k \O_{L_w}$ and $(\ell_L)_w(\pi_v)=[L_w: F_v]
\ell_v(\pi)$ so $(\ell_L)_w(\theta_w^{-1} \NN_w)=[L_w: F_v]
\ell_v(\theta_v^{-1} \NN_v)$. Thus,
\begin{equation*}
  \sum_{w|v} (\ell_L)_w(\theta_w^{-1} \NN_w) = m \ell_v(\theta_v^{-1}
  \NN_v)\;. 
\end{equation*}
If $w|v$ we have $\log_w = \log_v$ and $t_w = t_v\circ
\tr_{L_w/F_v}$. By Definition~\ref{logscext} we have
$\log_w(\theta_w(1))=\log_v(\theta_v(1))$ and thus clearly
\begin{equation*}
  \sum_{w|v} t_w(\log_w(\theta_w(1))) = m t_v(\log_v(\theta_v(1))\;.
\end{equation*}
This completes the proof.
\end{proof}
\begin{definition}
  Let $\NN$ be a line bundle over $\Spec(A)$ which is metrized over
  $F$. Then the Euler characteristic of $\NN$ is defined to be
  \begin{equation*}
    \chi(\NN)=\deg (\det_{\O_F} \NN)\;,
  \end{equation*}
  where $\det_{\O_F} \NN$ is metrized according to the log functions obtained
  from Lemma~\ref{goodlog}.
\end{definition}
The following two results are immediate consequences of Propositions
\ref{9.13}~and~\ref{9.14} respectively.
\begin{proposition}\label{chideg}
  For $\NN$ as above we have $\deg(\NN)= \chi(\NN)-\chi(A)$.
\end{proposition}
\begin{definition}\label{dualizing}
  Let $A$ and $L$ be as above. The dualizing module of $A$ over $\O_F$ is
  given by
  \begin{equation*}
    W_{A/\O_F}:=\{b\in L\;:\; \tr_{L/F} (b A)\subset \O_F\}\;,
  \end{equation*}
  metrized by the log functions induced from the inclusion into $L$.
\end{definition} 
\begin{proposition}\label{discform}
  We have $ \chi(A)= -\chi(W_{A/\O_F})=-\frac{1}{2} \deg(W_{A/\O_F}) $.
\end{proposition}
\begin{proposition}\label{intform}
  Let $D$ be an Arakelov divisor on $\XX$ and $E=\Spec(A)\subset \XX$ a
  horizontal divisor, with $A$ finite over $\O_F$. Then $D\cdot
  E=\deg(\O(D)|_E)$.
\end{proposition}
\begin{proof}
This is easily checked for an infinite fiber, so by linearity we may
assume that $D$ is an irreducible subscheme of codimension 1, and by a
moving lemma
on the generic fiber that the intersection of $D$ and $E$ with the
generic fiber have disjoint supports. It follows that  $\L:=\O(D)$ has a
global section $\O_{\XX} \xrightarrow{s} \L$, and this diagram serves as a
locally free resolution of $\O_D$. It also follows that $D$ and $E$
have proper intersection. Let $i:E \to \XX$ be the embedding and
$f:\XX \to \Spec(\O_F)$ the structure map.
We must compute the degree of the $\O_F$-line bundle
\begin{equation*}
  \M:= \det (\L|_E) \otimes (\det \O_E)^{-1} = \det ((\O_{\XX}
  \xrightarrow{s} \L)|_E)\;.
\end{equation*}
Here we implicitly must push down from $E$ to $\Spec(\O_F)$ along the
map which we can write as $f\circ i$.
The bundle $\M$ has a canonical section $s':\O_F \to \M$ induced by
the restriction to $E$ of the commutative diagram
\begin{equation*}
  \xymatrix{
   {\O_{\XX}} \ar[r]^{s} & \L\\
   {\O_{\XX}} \ar[r]^{\id}\ar[u]^{\id} & \O_{\XX} \ar[u]^{s}\\ 
  }
\end{equation*}
and the obvious triviality of the determinant of the bottom row.
Now we see that we can compute $\M$ as follows:
\begin{align*}
  \M &= \det \R f_\ast \R i_\ast (\O_{\XX} \to \L)|_E \\
  &= \det \R f_\ast ((\R i_\ast \O_E) \otimes (\O_{\XX} \to \L))
     \quad \text{by the projection formula}\\
  &= \det \R f_\ast (\shF^\bullet \otimes (\O_{\XX} \to \L))
\end{align*}
where $\shF^\bullet$ is a locally free resolution of $\O_E$. Since the
cohomology of $\shF^\bullet \otimes (\O_{\XX} \to \L)$ is an $\O_{\XX}$-module
supported exactly on the closed points of intersection between $D$ and
$E$, and since for such modules the map $f_\ast \to \R f_\ast$ is a
quasi-isomorphism, it follows that
\begin{equation*}
  \M \isom \bigotimes_{x\in D\cap E} (\otimes_i \det( f_\ast
  \operatorname{Tor}^{\O_{\XX}}_i(\O_D,\O_E))^{(-1)^i})\;.
\end{equation*}
Replacing $\O_{\XX} \to \L$ by $\O_{\XX} \to \O_{\XX}$ it is clear
that $s'$ is the alternating product of the maps induced by $0\to
\operatorname{Tor}^{\O_{\XX}}_i(\O_D,\O_E)$. For any place $v$
the determinant of the map $0\to \O_{F_v}/\pi_v^k \O_{F_v}$ tensored
with $F_v$ is such that the inverse image of $\det (\O_{F_v}/\pi_v^k
\O_{F_v})$ is $\pi_v^k \O_{F_v}$~\cite[Theorem~3 (vi)]{Knu-Mum76}. It
follows that the isomorphism
$\theta:=s'\otimes F: F\to \M\otimes F$ has
\begin{equation*}
  \theta_v^{-1}(\M_v) = \pi_v^k \O_{F_v},\quad k=\sum_{\substack{x\in
  D\cap E\\x\text{ above } v}} \sum_i (-1)^{i} \operatorname{length}
  (\operatorname{Tor}^{\O_{\XX}}_i(\O_D,\O_E))=\pair{D,E}_v\;.
\end{equation*}

Now we turn to the infinite contributions. Suppose $v|p$ and write
$E_v = \sum Q_j$. Then
\begin{equation*}
  \M_v=\det ((\O_{X_v} \xrightarrow{s} \L_v)|_{E_v}) = \otimes_j
  \det ((\O_{X_v} \xrightarrow{s} \L_v)|_{Q_j})\;.
\end{equation*}
The isomorphism we have chosen with $F_v$ is the tensor product of the
isomorphism of the j-th term with $F_v$, which is exactly the
isomorphism which was used in Definition~\ref{detlog}. By this
definition it is easy to see that 
\begin{equation*}
  \log (\theta_v(1))=\sum_j \tr_{F_v(Q_j)/F_v}
  \log_{\L_v}(s(1))(Q_j)=\sum_j \tr_{F_v(Q_j)/F_v} G_{D_v}(Q_j)=
  \pair{D,E}_v\;.
\end{equation*}
This completes the proof.

\end{proof}

\section{The adjunction formula and the Riemann-Roch theorem}
\label{sec:thms}

In this section we would like to show how some of the main theorems of
classical Arakelov theory have precise analogues in $p$-adic Arakelov
theory. In fact, after the work of the previous sections, the proofs
do not differ much from the proofs in the classical case. We have
chosen to follow the treatment of Lang~\cite{Lan88}.

We begin with the adjunction formula. Let $E\subset {\XX}$ be a horizontal
curve, with $E=\Spec(A)$ and $A$ finite over $\O_F$. Let
$\bfo_{E/\O_F}$ be the relative dualizing module. It is known that
\begin{equation}\label{geoadj}
  \bfo_{E/\O_F} = (\bfo_{{\XX}/\O_F}\otimes \O(E))|_E\;.
\end{equation}
Let $F(E)$ be the
function field of $E$. The residue map gives an injection
$\res:\Gamma(E,\bfo_{E/\O_F})\inject F(E) $.
\begin{theorem}[{\cite[Theorem~4.1, p. 94]{Lan88}}]
  The image of $\res$ is the dualizing module $W_{A/\O_F}$ of
  Definition~\ref{dualizing}.
\end{theorem}
Note that this dualizing module is taken in this definition without
its metric. This metric figures in the next definition.
\begin{definition}\label{paddisc}
  The ($p$-adic) \emph{discriminant} of $E$ is $d(E)=\deg(W_{A/\O_F})$.
\end{definition}
This definition does not take into account the embedding of
$E$ in ${\XX}$. Let now $v$ be a prime above $p$ and let $\fvp$ be an
algebraic closure of $F_v$. By definition we have a Green function
$G_v$ on $X\otimes \fvp$. Over $\fvp$ the divisor $E$ splits as a sum
of distinct points $E\otimes \fvp = \sum_{j=1}^e P_j$.
\begin{definition}\label{infdisc}
  The \emph{discriminant above $v$} of $E$ in $X$ is 
  \begin{equation*}
    d_v(E,\XX)=\sum_{i\ne j} G_v(P_i,P_j)\;.
  \end{equation*}
  The \emph{infinite discriminant} is defined as
  \begin{equation*}
    d_\infty(E,\XX)=\sum_{v|p} d_v(E,X)\;.
  \end{equation*}
\end{definition}

The adjunction formula is now the following statement.
\begin{theorem}\label{adjunction}
  Let $E$ be a horizontal divisor on $X$. Then
  \begin{equation*}
    \bfo_{X/\O_F} \cdot E + E\cdot E = d(E)+d_\infty(E,\XX)\;.
  \end{equation*}
\end{theorem}
\begin{proof}
From \eqref{geoadj} and the fact that both $\bfo_{X/\O_F}$ and $\O(E)$
have natural metrics, we obtain a metric on $\bfo_{E/\O_F}$. With
respect to this metric it follows from Proposition~\ref{intform} that
  \begin{equation*}
    \bfo_{X/\O_F} \cdot E + E\cdot E = \deg(\bfo_{E/\O_F})\;.
  \end{equation*}
As we defined them, $\bfo_{E/\O_F}$ and $W_{A/\O_F}$ are the same
module but with different metrics. The difference in their degree is
thus the sum of the differences between their log functions. Consider
a place $v|p$ of $F$ and a point $P_i$ in $E\otimes \fvp$ as before.
The log function on the fiber at $P_i$ of $W_{A/\O_F}$ is such that
the residue map to $\fvp$ is an
isometry. On the other hand, the fiber at the same
point of $\bfo_{E/\O_F}$ is viewed as the fiber of
$\bfo_{\xvb/\fvp}\otimes \O(\sum P_j)$. The log function on the fiber
of $\bfo_{\xvb/\fvp}\otimes \O(P_i)$ is such that the residue is an
isometry and the points $P_j$ for $j\ne i$ contribute an added term of
$G_v(P_i,P_j)$. The result is now clear.
\end{proof}

Suppose now that $\XX$ is an arithmetic surface over $\O_F$ and that
$\L$ is a metrized line bundle over $\XX$. Let 
\begin{equation*}
  \M= \lambda(\L):= \det H^0(\XX,\L) \otimes (\det H^1(\XX,\L))^{-1}\;.
\end{equation*}
Then,
for any place $v|p$ of $F$ we have $\M_v = \lambda(\L_v)$ and by
Proposition~\ref{detcond}
was done before it acquires a log function. Thus, $\M$ is a metrized
line bundle on $\O_F$ and we can define
\begin{equation*}
  \chi(\L)=\deg (\lambda(\L))\;.
\end{equation*}
The following lemma is the $p$-adic analogue of a well known result in
classical Arakelov theory, and is an immediate consequence of the
multiplicativity of the determinant and the behavior of the Faltings
volume with respect to adding points given in part \eqref{detcond3} of
Proposition~\ref{detcond}.
\begin{lemma}\label{9.5}
  If $D$ is an Arakelov divisor and $E$ is a horizontal divisor, then
  \begin{equation*}
    \chi(\O(D+E))=\chi(\O(D))+\chi(\O(D+E)|_E)-\frac{1}{2}d_\infty(E,\XX)\;.
  \end{equation*}
\end{lemma}
\begin{theorem}\label{RR}
  We have the following Riemann-Roch formula:
  \begin{equation*}
    \chi(\L) - \chi(\O_{\XX})=\frac{1}{2} \L \cdot (\L-\omega_{\XX})\;. 
  \end{equation*}
\end{theorem}
\begin{proof}
We follow the proof given by Lang. One can assume that $\L$ is of the
form $\O(D)$ for some Arakelov divisor $D$. Then one checks that the
validity of the theorem is unchanged if one adds or subtracts from $D$
a divisor. The two cases of a fiber at infinity and of a vertical
divisor are essentially the same as in the classical case so we leave
them for the reader. The case of adding a horizontal divisor is
treated exactly as in Lang. We reproduce the proof to see that we have
all the ingredients (with slightly different notation). We have
\begin{align*}
  \chi(\O(D)|_E)&= \deg(\O(D)|_E) + \chi(\O_E)\quad \text{by
  Proposition~\ref{chideg}}\\
  &=D \cdot E  + \chi(\O_E)\quad \text{by
  Proposition~\ref{intform}}\\
  &= D\cdot E - \frac{1}{2} d(E)\quad \text{by
  Proposition~\ref{discform} and Definition~\ref{paddisc}}\\
  &= D\cdot E -\frac{1}{2} \left(E\cdot E + \bfo_{\XX/\O_F}\cdot E -
  d_\infty(E,\XX) \right)
\end{align*}
by the adjunction formula (Theorem~\ref{adjunction}). If we replace
$\O(D)$ by $\O(D+E)$ the left hand side of the Riemann-Roch formula
changes by
\begin{align*}
   \chi(\O(D+E))&-\chi(\O(D))=
   \chi(\O(D+E)|_E)-\frac{1}{2}d_\infty(E,\XX)\quad
   \text{by Lemma~\ref{9.5}}\\
   &=(D+E) \cdot E -\frac{1}{2} \left(E\cdot E + \bfo_{\XX/\O_F}\cdot E -
  d_\infty(E,\XX) \right))-\frac{1}{2}d_\infty(E,\XX)\\
  &=D\cdot E +\frac{1}{2} E\cdot (E- \bfo_{\XX/\O_F})\;,
\end{align*}
which is exactly the amount by which the right hand side changes.
\end{proof}
\begin{remark}
It is evident from the proof of the Riemann-Roch theorem that it is
too simple to depend on a particular normalization of the Green
function $G$ or of the log function on the determinant of
cohomology. This latter independence is clear. Here we would like to
check the independence of the Green function directly, since this 
requires keeping careful track of all normalizations. Suppose then
that we have two Green functions $G_1$ and $G_2=G_1+1$ at the place
$v$ and the Green functions at the other places are the same. It suffices to
consider this case since all contributions will be linear in the
constant $G_2-G_1$. We check how the two sides of Theorem~\ref{RR} change when
we make this change, beginning with the right hand side. With respect
to these two functions we have the 
following quantities that change: The intersection
pairing
$\pair{~,~}_i$, the Arakelov Chern class $c_i$ and the canonical
class $\bfo_i$, for $i=1,2$. Let $d=$ be the degree of $\L$ on the
generic fiber.

The relation between the intersection products is that 
\begin{equation*}
  \pair{D,E}_2= \pair{D,E}_1+ \deg(D_F)\cdot \deg(E_F)\;.
\end{equation*}
Since $G_2(D,\bullet)=G_1(D,\bullet)+\deg D_F$ it follows from
Definition~\ref{ilogl} that
the relation between the $\ilog$ characters at the place $v$ is
$\ilog_{\L,2}=\ilog_{\L,1}-\deg(\L)$ and therefore
\begin{equation*}
  c_2(\L)=c_1(\L)- dX_v\;.
\end{equation*}
Finally, Proposition~\ref{charofgr} implies that
\begin{equation*}
  \log_{{\bfo_2}_v} = \log_{{\bfo_2}_v}-1\;
\end{equation*}
from which it follows that
\begin{equation*}
  c_2(\bfo_2)=c_1(\bfo_1)- (2g-2+1) X_v
\end{equation*}
Thus, (twice) the right hand side of the Riemann-Roch formula changes
as follows:
\begin{align*}
  &\phantom{=} \pair{c_2(\L),c_2(\L)-c_2(\bfo_2)}_2\\&=
  \pair{c_1(\L)-dX_v,c_1(\L)-c_1(\bfo_1)-(d-2g+1)X_v}_1+
  d(d-2g+2)\\ &=
  \pair{c_1(\L),c_1(\L)-c_1(\bfo_1)}_1-d(d-2g+2)-d(d-2g+1)
  +d(d-2g+2)\\ &=\pair{c_1(\L),c_1(\L)-c_1(\bfo_1)}_1-d(d-2g+1)\;.
\end{align*}

Now we turn to the left hand side. We have
\begin{equation*}
  \chi(\L)- \chi(\O)=(\chi(\L)-\chi(\O(D)))+(\chi(\O(D))-\chi(\O))
\end{equation*}
where $D$ is the finite part of the Chern class of $\L$.
The first summand reflects the different metric between $\O(D)$ and
$\L$. The change of $G$ adds $d$ to $\log_{\O(D)}$, which, in view of
\eqref{detcond2} of Proposition~\ref{detcond}, subtracts $(d+1-g)d$
from the first summand. On the other hand,
it follows from \eqref{detcond3} of Proposition~\ref{detcond} that the
second summand gets $1+2+\cdots d= d(d+1)/2$ added. So overall the
left hand side is reduced by
\begin{equation*}
  (d+1-g)d  - \frac{d(d+1)}{2} = d \left( \frac{d(d+1)}{2}-g\right)\;,
\end{equation*}
which is exactly  what gets subtracted from the right hand side.
\end{remark}

\appendix

\section{The universal vectorial extension of a Jacobian}
\label{sec:univ}

In this appendix we review the theory of the universal vectorial
extension of the Jacobian of a curve and prove several algebraic
results that will be required in the main text. Probably, everything
is well known but we do not know of a reference. The general theory of
vectorial
extensions of abelian varieties is to be found in~\cite{Maz-Mes74}. It
is utilized for Jacobians by Coleman in~\cite{Col90,Col91}, but our
treatment is independent of his.

Let $C$ be a curve over a base scheme $S$. Let $K(C)^\times$ be the
sheaf
\begin{equation*}
  \bigoplus_{\operatorname{cod} x=0} i_x K(x)^\times
\end{equation*}
where $i_x$ is the embedding of $x$ in $C$ and $k(x)^\times$ is the
 multiplicative group of the residue field of $x$.
 We let $\tknd{C/S}$ be the complex of sheaves
\begin{equation*}
  \O_C^\times \xrightarrow{\dd} \Omega_{C/S}^1 \oplus K(C)^\times\;,\;
  \dd(f)=(\dlog(f),f)\;.
\end{equation*}
\begin{definition}
  The space of differentials of the third kind on $C$ relative to $S$
  is the group $H^1(C,\tknd{C/S})$.
\end{definition}
To see why this definition captures differentials of the third kind we
compute this cohomology with the help of a Zariski covering $\{U_i\}$
of $C$. A one cocycle is given by
\begin{equation*}
  (g_i\in K(C)^\times(U_i),\; \omega_i \in \Omega_{C/S}(U_i),\;
  f_{ij} \in \O(U_{ij}))
\end{equation*}
such that
\begin{equation*}
\omega_i-\omega_j=
  \dlog(f_{ij})\text { and } f_{ij}=\frac{g_i}{g_j}
\end{equation*}
Given such a cocycle, we recover a form of the third kind by taking
$\omega_i - \dlog(g_i)$ on $U_i$ and noticing that the conditions
guarantee that these glue together. The resulting form has by
definition logarithmic singularities. More conceptually, when $C$ and
$S$ are spectra of fields we have an isomorphism $\tknd{C/S}\to
\Omega_{C/S}^1[1]$ given by $(\omega,g)\mapsto
\omega- \dlog(g)$. Thus we obtain from a form of the third kind a
differential at the generic point by first restricting and then
applying this isomorphism on cohomology.

In~\cite[I.3.1.7]{Maz-Mes74} the multiplicative de Rham complex of $C/S$ is
defined to be the complex $\Omega_{C/S}^\times= (\O_C^\times
\xrightarrow{\dlog} \Omega_{C/S}^1)$ (in loc.\ cit.\ it is extended
further to the right, which we do not have to do). There is an obvious
short exact sequence
\begin{equation*}
  0\to K(C)^\times[1] \to \tknd{C/S} \to \Omega_{C/S}^\times \to 0\;.
\end{equation*}
Taking cohomology we obtain the short exact sequence
\begin{equation*}
  K(C)^\times \to H^1(C,\tknd{C/S}) \to H^1(\Omega_{C/S}^\times) \to 0\;.
\end{equation*}
It is known that when sheafifying the right term of the above sequence
one obtains a functor represented by the universal vectorial extension
$G_X$ of the Jacobian $J$ of $C$. The map on the left sends
a rational function $f$ to the form of the third kind $-\dlog(f)$.

Recall from Section~\ref{sec:double} the definition of differentiation
of differential forms with respect to a vector field. We are going to
refine this to a differentiation from a family of forms of the third
kind to a family of forms of the second kind, a notion defined as follows.
\begin{definition}
  A family of forms of the second kind on $C/S$ is an element of
  $H^1(C,\sknd{C/S})$, where $\sknd{C/S}$ is the complex
  \begin{equation*}
    \O_C \xrightarrow{\dd} \Omega_{C/S}^1 \oplus K(C)\;,\; \dd(f)=(\dd
    f, f)\;.
  \end{equation*}
\end{definition}
We have an obvious short exact sequence
\begin{equation*}
  0\to K(C)[1] \to \sknd{C/S} \to \Omega_{C/S}^\bullet \to 0\;.
\end{equation*}
Taking cohomology we obtain
\begin{equation*}
  K(C) \to H^1(C,\sknd{C/S}) \to H^1(C,\Omega_{C/S}^\bullet) \to 0\;.
\end{equation*}
When $S=\Spec(K)$ this map describes the representation of the first
de Rham cohomology of $C$ as the quotient of the space of forms of the
second kind by the differentials of rational functions.
\begin{definition}
  Let $\partial/\partial t$ be a vector field on $C$. We define a map
  \begin{equation*}
    \frac{\partial}{\partial t}: \tknd{C/S} \to \sknd{C/S}\;,
  \end{equation*}
  given in degree $0$ by
  \begin{equation*}
    f\mapsto \frac{\frac{\partial f}{\partial t}}{f}
  \end{equation*}
  and in degree $1$ by
  \begin{equation*}
    (\omega,f)\mapsto \left(\frac{\partial}{\partial t} \omega,
    \frac{\frac{\partial f}{\partial t}}{f}\right)\;.
  \end{equation*}
\end{definition}
It is easy to check that this is indeed a map of complexes. On $H^1$
it gives a map, which we continue to call $\partial/\partial t$, from
forms of the third kind to forms of the second kind. It is further
easy to check that viewing both forms of the third and second kind as
differential forms on the generic point, this map is just
differentiation of forms with respect to the restriction of
$\partial/\partial t$ to this point.

\begin{proposition}
  Consider the family $C[\varepsilon]/S[\varepsilon]$ where
  $\varepsilon^2=0$. Then we have the following commutative
  diagram
    \begin{equation*}
    \xymatrix{
    {\Ker \left(H^1(C,\tknd{C[\varepsilon]/S[\varepsilon]})\to
    H^1(C,\tknd{C/S})\right)} \ar[d] \ar[r]^-{\partial/\partial \varepsilon}
    & H^1(C,\sknd{C/S})\ar[d] \\
    {\Ker \left(H^1(C,\Omega_{C[\varepsilon]/S[\varepsilon]}^\times)\to
    H^1(C,\Omega_{C/S}^\times)\right)} \ar[r]^-{\sim} & \hdr^1(C/S)
          }
  \end{equation*}
  In this diagram the top horizontal map is differentiation composed with
  restriction to $C$ and the vertical map on the right sends a
  differential of the third kind to its de Rham cohomology class.
\end{proposition}
\begin{proof}
Suppose we
have a form of the third kind on $C[\varepsilon]/S[\varepsilon]$ whose
restriction to $C/S$ is $0$. The
element $\varepsilon$ provides a canonical vector field on
$C[\varepsilon]$ and we would like to compute the derivative of this
form of the third kind with respect to $\varepsilon$. First we notice
that there is a commutative diagram with split short exact sequence
rows of sheaves on $C$ 
(compare the proof of Proposition I.4.1.4 in~\cite{Maz-Mes74})
\begin{equation*}
  \xymatrix{
    0\ar[r]& \Omega_{C/S}^\bullet \ar[r]\ar@2{-}[d] &
    \tknd{C[\varepsilon]/S[\varepsilon]} 
       \ar[r] \ar[d]&\tknd{C/S} \ar[r]\ar[d]& 0\\
    0\ar[r]& \Omega_{C/S}^\bullet \ar[r]&
    \Omega_{C[\varepsilon]/S[\varepsilon]}^\times 
       \ar[r] &\Omega_{C/S}^\times \ar[r]& 0\\
  }
\end{equation*}
In this diagram, the vertical maps are the ones defined before. The
top left horizontal map in degree $0$ sends $f$ to $1+\varepsilon f$
and in degree $1$ sends $\omega$ to $ \varepsilon \omega$. It is now
easy to see that the composed map
\begin{equation*}
  \Omega_{C/S}^\bullet \to \tknd{C[\varepsilon]/S[\varepsilon]}
  \xrightarrow{\partial/\partial \varepsilon}
  \sknd{C[\varepsilon]/S[\varepsilon]} \to
  \sknd{C/S} \to \Omega_{C/S}^\bullet
\end{equation*}
is the identity map. Indeed, in degree $0$ it first sends $f$ to
$1+\varepsilon f$. Then log differentiating with respect to $\varepsilon$
sends this to $f/(1+\varepsilon f)$ and this is then sent back to $f$
by the map that kills $\varepsilon$ and this is sent to $f$ again. In
degree $1$ a form $\omega$ is sent to $(\varepsilon \omega,0)$,
differentiation with respect to $\varepsilon$ sends this to
$(\omega,0)$, then to $(\omega,0)$ again and finally to $\omega$. By
taking cohomology we obtain the result.
\end{proof}
\begin{corollary}\label{A5}
  Let $C/K$ be a complete curve, $T/K$ a variety, $0\in T(K)$ a fixed
  point. Let $G_C$ be the universal vectorial extension of $J(C)$. Let
  $(\eta_t)_{t\in T}$ be a family of forms of the third
  kind on $C\times T/T$ and let $\indi:T\to G_C$
  be the induced map. Suppose $\eta_0=0$, which implies
  $\indi(0)=0$. Let $\partial/\partial t$ be a vector field on
  $T$. Then, the form of the second kind on $C$, $(\partial
  \eta_t/\partial t)|_{t=0}$, represents the de Rham cohomology class
  \begin{equation*}
    (\dd \indi) (\frac{\partial}{\partial t}|_{t=0}) \in
    \operatorname{Lie}(G_C)\isom \hdr^1(C/K)\;.
  \end{equation*}
\end{corollary}
\begin{proof}
  This follows from the previous proposition by restricting to the
  infinitesimal neighborhood of $0$ in $T$ and interpreting the result.
\end{proof}

Finally, to use the previous corollary, we want to show that we can at
least locally lift elements of the universal vectorial extension to
forms of the third kind.
\begin{lemma}\label{loclift}
  Consider $C$, $T$ and $0$ as in the corollary. Let $\indi: T\to G_C$
  be a map. Then there exist a neighborhood $U$ of $0$ in $T$ and a
  family of forms of the third kind  $(\eta_t)_{t\in U}$ inducing $\indi_U$.
\end{lemma}
\begin{proof}
By~\cite{Maz-Mes74} the map $\indi$ is locally induced by a line
bundle $\L$ on
$T\times C$ together with a relative connection $\nabla$ on it. To
obtain a family of forms of the third kind one takes a section $s$ of
$\L$ and compute the form of the third kind $\nabla(s)/s$. We may
choose $U$ and the section in such a way that $s$ is invertible on
$U$, hence the family of forms is defined on $U$.
\end{proof}
\section{Relations with the theory of Colmez}
\label{sec:colmez}

In~\cite{Colm96} Colmez developed a theory of $p$-adic integration
using Abelian varieties. In this theory there is also a notion of
Green functions. The purpose of this section is to compare this notion
of Green functions with the one we have been developing here. Note
that Colmez is working over $\C_p$ while we are working over $\qpc$.

Let $A$ be an abelian variety over $\qpc$. For each non-negative integer
$n$ Colmez defines a correspondence $\DD^n :A^{n+1} \to A$ (a kind of
difference operator) as follows:
Define, for $I\subset \{1,\ldots,n\}$, $m_I: A^{n+1}\to A$ by
\begin{equation*}
  m_I(x,h_1,\ldots,h_n)= x+\sum_{i\in I} h_i  
\end{equation*}
Then,
\begin{equation*}
  \DD^n := \sum_{I\subset \{1,\ldots,n\}} (-1)^{n-|I|} m_I^\ast\;.
\end{equation*}
The following easy lemma, taken from~\cite{Colm96}, gives an
alternative recursive description of
$\DD^n$.
\begin{lemma}
  Let
  \begin{align*}
    \pi_n, m_n: A^{n+1}\to A^n,\quad \pi_n (x,h_1,\ldots, h_n)&=
    (x,h_1,\ldots, h_{n-1}),\\ m_n (x,h_1,\ldots, h_n)&=
    (x+h_n,h_1,\ldots, h_{n-1})\;.
  \end{align*}
  Then, $\DD^n$ is given recursively by the formulas
  \begin{equation*}
    \DD^0 = \id,\quad \DD^n = (\pi_n^\ast - m_n^\ast)\circ \DD^{n-1} \;.
  \end{equation*}
\end{lemma}
It is immediate to see that for any $i\in \{1,\ldots,n\}$ restriction
to $h_i=0$ composed with $\DD^n$ equals $0$.

Now let $\L$ be a line bundle on $A$. The Theorem of the cube implies
the existence of an isomorphism 
\begin{equation}
  \label{eq:cube}
  \DD^3 \L \isom \O_{A^4}\;.
\end{equation}
We normalize this isomorphism by requiring that it restricts to the
identity isomorphism on each $\{h_i=0\}$.
\begin{proposition}
  For any log function $\llog$ on $\L$ the isomorphism \eqref{eq:cube}
  is an isometry of metrized line bundles.
\end{proposition}
\begin{proof}
Since $\DD^n$ is a difference operator it follows easily that $\DD^n$
kills $\hdr^1(A)$ for $n\ge 2$ and $\DD^n$ kills $\Ht(A)$ for $n\ge
3$. In particular, the curvature of $\DD^3 \L$ is $0$. It follows that
the differentials of the log functions on the two sides of
\eqref{eq:cube} can differ by at most a holomorphic differential on
$A^4$. This has the form $\sum_{i=0}^3 p_i^\ast \omega_i$ with $\omega_i \in
\Omega^1(A)$ and $p_i:A^4\to A$ the projection on the $i$th
coordinate. But the restriction of this form to $\{h_i=0\}$ is $0$,
showing that $\omega_i=0$. Thus, \eqref{eq:cube} is an isometry up to
scaling and again restricting to $\{h_i=0\}$ shows that it is in fact
an isometry.
\end{proof}

We now compare this result with Proposition~I.2.8
of~\cite{Colm96}. Let $s$ be a section of the line bundle $\L$ and let
$D$ be the divisor of $s$. For any log function $\llog$ on $\L$ let
$G_D=\llog(s)$. The Theorem of the cube implies that $\DD^3 D$ is a
principal divisor and Colmez chooses a rational function $f_D^{(4)}$
normalized in such a way that its restriction to $\{h_i=0\}$ is
$1$. This is clearly just the image of $\DD^3 s$ under the canonical
choice of~\eqref{eq:cube}. It is now immediate that $\DD^3 G_D = \log
f_D^{(4)}$, which is the defining property of the Green function of
the divisor $D$ in Colmez's definition. By the properties of Coleman
integration it is easily seen that $G_D$ is locally analytic outside
$D$ and has logarithmic singularities along $D$. It is therefore the
Green function of Colmez. The kernel of the cup product map $\cup:
\Ht(A) \to \hdr^2(A)$ is exactly $\symm^2 \Omega^1(A)$. Thus,
different choices for $\llog$, and consequently for $G_D$, differ by
the constant of integration,
by the integral of a holomorphic form on $A$, and by integrals
corresponding to elements of $\symm^2 \Omega^1(A)$, i.e., integrals of
the form
\begin{equation*}
  \int (\omega \int \omega) = \frac{1}{2} (\int \omega)\cdot (\int \omega)\;.
\end{equation*}
In other words, $G_D$ is unique up to a polynomial of degree $2$ in
the integrals of holomorphic forms on $A$, which are the logarithms of
$A$ in Colmez's terminology, and this is exactly the indeterminacy in
Colmez's Green functions. To sum up, we have proved
\begin{proposition}
  Let $\L$ be a line bundle on $A$, $s$ a section of $\L$ and $D$ the
  divisor of $s$. The collection of Green functions for $D$ defined by
  Colmez is the same as the collection of functions $\llog(s)$ for all
  possible log functions $\llog$ on $\L$. In particular, the Green
  functions of Colmez are Coleman functions.
\end{proposition}

\end{document}